\documentclass{amsart}
\usepackage{amscd}
\usepackage{amssymb}
\usepackage{hyperref}
\usepackage{cleveref}
\usepackage{url}
\usepackage{stmaryrd}

\RequirePackage{tikz-cd}
\RequirePackage{amssymb}
\usetikzlibrary{calc}
\usetikzlibrary{decorations.pathmorphing}

\RequirePackage{tikz}

\newcommand{\comment}[1]{}	

\newcommand{\A}{\mathbb{A}}

\input xypic
\xyoption{all} 
\numberwithin{equation}{section}

\theoremstyle{plain}
\newtheorem{theorem}{Theorem}[subsection]
\newtheorem{corollary}[theorem]{Corollary}
\newtheorem{lemma}[theorem]{Lemma}
\newtheorem{proposition}[theorem]{Proposition}

\theoremstyle{definition}
\newtheorem{definition}[theorem]{Definition}
\newtheorem{example}[theorem]{Example}
\newtheorem{observation}[theorem]{Observation}
\newtheorem{construction}[theorem]{Construction}
\theoremstyle{remark}
\newtheorem{remark}[theorem]{Remark}

\begin{document}
\bibliographystyle{plain}

\title{Resolution of embedded toric $\Lambda$-schemes}

\author{K. Machida}
\thanks{Machida was supported by the David Lachlan Hay Memorial Fund through the University of Melbourne Faculty of Science Postgraduate Writing-Up Award.}
\address{School of Mathematics and Statistics, University of Melbourne,
VIC 3010, Australia}
\email{machida@student.unimelb.edu.au}

\date{\today}

\begin{abstract}
A key example in Borger's theory of $\Lambda$-structure is toric $\Lambda$-structure. We prove a resolution of singularities result for embedded toric $\Lambda$-schemes by applying an algorithm of Bierstone and Milman for toric varieties over perfect fields. This paper is based on work from the author's PhD thesis.
\end{abstract}

\maketitle

\section*{}

In this paper we prove a resolution theorem for embedded toric $\Lambda$-schemes analogous to the embedded desingularization of toric varieties by Bierstone and Milman (see \cite[Theorem 1.1]{BM1} and \cite[Theorem 14.1]{monoid}). Toric $\Lambda$-schemes are the $\Lambda$-schemes (see \cite{borger2009lambdarings}) with toric $\Lambda$-structure arising from the $\mathbb{Z}$-realization of pointed monoid schemes as defined by Corti\~nas, Haesemeyer, Walker and Weibel in \cite{monoid}. Our method involves defining analogous objects in the context of $\Lambda$-schemes and showing that certain properties commute with fibers. This reduces the resolution to order reduction of $\Lambda$-marked monomial ideals developed in \cite{OrderReductionLambdamarkedMonomialIdeals}. In particular, we show that the two invariants, order and Hilbert-Samuel functions, that are used to stratify the singular locus behave well with respect to fibers. Moreover, these stratifications are shown to be $\Lambda$-schemes themselves. We may interpret this result, in the often more philosophical than mathematical spirit of $\mathbb{F}_1$-geometry as evidence that the singular points we blow up really should be viewed as coming from $\mathbb{F}_1$. Finally, the results of the paper also show that the $\Lambda$-resolution can be viewed as a family of resolutions over $\textrm{Spec} \,\mathbb{Z}$, with each fiber the resolution achieved by applying \cite[Theorem 1.1]{BM1} directly. 
\\In \cref{C7_3} and \cref{C7_4} we describe the main $\Lambda$-schemes we consider. These are analogues of toric varieties and we show that there are torus actions and invariant closed $\Lambda$-subschemes coming from the underlying monoid schemes. The closed embedding of toric $\Lambda$-schemes is described in \cref{C7_5} as locally given by binomial equations. Section \ref{C7_6} describes blowing up and taking strict transforms from the perspective of local binomial equations and underlying monoid schemes. In \cref{C7_8} we use Hasse derivatives and properties of cones and fans to prove the main theorem in the simpler setting of a hypersurface. The arguments in this section follow \cite[Section 5]{BM1} though we use Hasse derivatives to give stratifications with the structure of a $\Lambda$-scheme. Section \ref{C7_9} shows how order at points in toric $\Lambda$-schemes behave with respect to fibers and \cref{C7_11} uses the standard basis elements of \cref{C7_9} to do the same with Hilbert-Samuel functions. Together \cref{C7_9} and \cref{C7_11} prove the necessary properties that allow us to use the inductive arguments of Bierstone and Milman. Finally, in \cref{C7_12} we prove the main theorem using the arguments of Bierstone and Milman applied to each fiber and order reduction for $\Lambda$-marked monomial ideals.


\begin{theorem} 
    \label[theorem]{theorem7201} 
    Let $X_{MSch}$ be a connected separated cancellative torsionfree monoid scheme of finite type, embedded as an unpointed closed subscheme in a smooth toric monoid scheme $M_{MSch}$. Let $X \hookrightarrow M$ be the $\Lambda$-equivariant embedding of toric $\Lambda$-schemes given by the $\mathbb{Z}$-realization of the monoid scheme embedding, where $M$ is smooth over $\textrm{Spec} \, \mathbb{Z}$. Then there is a finite sequence of $\Lambda$-equivariant blow-ups of $M$,
    $$
    M=M_0 \xleftarrow[]{\pi_1} M_1 \leftarrow \dots \xleftarrow[]{\pi_{l+1}} M_{t+1},
    $$
    such that:
    \begin{enumerate}
        \item The centre $D_j$ of each blowing-up $\pi_{j+1}$ is a smooth $T_{M_j}$-invariant $\Lambda$-subscheme of $M_j$.
        \item Set $X_0=X$. For each $j=0, \dots t$, let $X_{j+1}$ denote the strict transform of $X_j$ by $\pi_{j+1}$. Then each $C_j := D_j \cap X_j$ is a $T_{X_j}$-invariant $\Lambda$-subscheme of the $\mathbb{Z}$-flat $\Lambda$-scheme $X_j$. Moreover, $C_j$ is smooth over $\textrm{Spec} \, \mathbb{Z}$.
        \item For each $j=0, \dots t$, $C_j \subset Sing_{\mathbb{Z}}(X_j)$, where $Sing_{\mathbb{Z}}(X_j)$ is the locus of points that are not smooth in their fibers. Moreover, $X_{t+1}$ is smooth over $\textrm{Spec} \, \mathbb{Z}$. 
    \end{enumerate}
\end{theorem}

\begin{remark}
    The theorem is an adaptation of a resolution theorem for embedded toric varieties over a perfect field by Bierstone and Milman (see \cite[Theorem 1.1]{BM1}). The conditions for the embedding of monoid schemes are those of \cite[Theorem 14.1]{monoid} with additional assumptions (connected and unpointed). As a result of the monoid scheme structure Theorem \ref{theorem7201} gives a $\Lambda$-resolution that addresses the limitations of \cite[Observation 5.0.6]{OrderReductionLambdamarkedMonomialIdeals}.
\end{remark}


In fact, a simple observation and the proof of \cite[Theorem 14.1]{monoid} would be enough to prove the theorem. \cite[Theorem 14.1]{monoid} is argued by applying Bierstone and Milman’s theorem to the $k$-realization ($k$ a field) of the monoid scheme embedding to produce a sequence of monoid scheme blow-ups that results in a smooth monoid scheme. Applying $\mathbb{Z}$-realization to this monoid scheme blow-up sequence would give our $\Lambda$-equivariant blow-up sequence resulting in a smooth $\Lambda$-scheme. However, such an observation would not tell us how the data of each $\Lambda$-scheme embedding is used to choose the centres of blowing up and would not give any insight into how we might generalise to other types of $\Lambda$-schemes. Further, in order to work with $\mathbb{Z}$-realizations of these monoid schemes as in Bierstone and Milman we require slightly stronger conditions than just those in \cite[Theorem 14.1]{monoid}. 
\\
\\\textbf{Terminology.}
\\
\\\textbf{Schemes.} A scheme will be a ringed space $(X, \mathcal{O}_X)$ that can be covered by affine schemes.
\\
\\\textbf{Varieties.} An algebraic variety over a field $k$ will be a separated reduced scheme of finite type over $k$.
\\
\\\textbf{$\mathbb{Z}$-torsion free and $\mathbb{Z}$-flatness.} Let $R$ be a ring. $R$ will be $\mathbb{Z}$-torsion free if it is so as a $\mathbb{Z}$-module, i.e. as an abelian group under addition. Let $\textrm{Spec} \,R$ be an affine scheme, then as $\mathbb{Z}$ is an integral domain $\textrm{Spec} \,R$ being flat over $\mathbb{Z}$ is equivalent to $R$ being $\mathbb{Z}$-torsion free. For affine schemes, we will often use these terms interchangeably, emphasising the structure of the ring or the scheme depending on the context.
\\
\\\textbf{Monoids.} All monoids in this paper will be commutative.


\section{Toric \texorpdfstring{$\Lambda$}{Λ}-schemes}
\label{C7_3} 
We describe the $\Lambda$-scheme analogues of toric varieties and actions of tori. A \textit{toric monoid scheme} is a separated connected torsionfree monoid scheme of finite type (see \cite{monoid}). Our definition differs from \cite[Section 4]{monoid} in that our toric monoid schemes need not be normal. We show a toric monoid scheme has an action of a torus, justifying the use of the term toric. A \textit{normal toric monoid scheme} is then a separated connected torsionfree normal monoid scheme of finite type, which is the definition of a toric monoid scheme in \cite[Section 4]{monoid}.

\begin{definition}
    We will call a scheme $X$ such that $X \cong (X_{MSch})_{\mathbb{Z}}$ is the $\mathbb{Z}$-realization of a a toric monoid scheme $X_{MSch}$ a \textit{toric $\Lambda$-scheme} (see \cite[Section 14]{monoid}). We describe an action of a torus on $X$ that justifies the use of the term toric.
\end{definition}

\begin{observation}
    By definition $X$ is covered by affine open subschemes $\textrm{Spec} \, \mathbb{Z}[A]$ for $A$ cancellative torsionfree monoids of finite type (see \cite[Section 1]{monoid}). Then as $\mathbb{Z}$ is an integral domain and the group completion $A^{+}$ is torsionfree, $\mathbb{Z}[A]$ is an integral domain and $X$ an integral scheme. Further, so is each fiber over $(p) \in \textrm{Spec}  \, \mathbb{Z}$, $X_p \cong (X_{MSch})_{\mathbb{F}_p}$ (see \cite[Lemma 5.10]{monoid}).
\end{observation}

\begin{lemma}
    \label[lemma]{lemma103} 
    Let $X \cong (X_{MSch})_{\mathbb{Z}}$ for $X_{MSch}$ a separated cancellative monoid scheme of finite type. Then $X$ is smooth over $\textrm{Spec} \, \mathbb{Z}$ if and only if $X_{MSch}$ is a smooth monoid scheme.
\end{lemma}

\textit{Proof.} Suppose $X_{MSch}$ is a smooth monoid scheme, i.e. its stalks are the smash product of a free group of finite rank and a free monoid of finite rank. Then $X$ can be covered by affine open subschemes of the form $\textrm{Spec} \, \mathbb{Z}[x_1, \dots, x_r, y_{r+1}^{\pm}, \dots, y_n^{\pm}]$ and hence $X$ is smooth over $\textrm{Spec} \, \mathbb{Z}$ (see \cite[Theorem 25.2.2]{FOAG}). Conversely, suppose $X$ is smooth over $\mathbb{Z}$, then each fiber $X_p \cong (X_{MSch})_{\mathbb{F}_p}$ is smooth over $\mathbb{F}_p$ hence by \cite[Proposition 6.5]{monoid}, $X_{MSch}$ is smooth. \qedsymbol


\subsection{Torus actions}
\begin{construction}
    \label{construction7311} 
    Let $X$ be a toric $\Lambda$-scheme corresponding to $X_{MSch}$. Since $X_{MSch}$ is cancellative and connected it is the closure of a unique minimal point $\eta$ (see \cite[lemma 2.3]{monoid}), which is contained in every open affine chart $\textrm{MSpec}(A)$ of $X_{MSch}$ as the prime ideal $0$ in the cancellative finitely generated mononid $A$. The localization at the prime ideal corresponding to $\eta$, which we denote $A_{\eta}$ corresponds to a distinguished open subscheme $\textrm{MSpec}(A_{\eta})$ of $\textrm{MSpec}(A)$ (\cite[Lemma 1.5]{monoid}) given by localization at the multiplicative set consisting of all generators of $A$. Thus $\textrm{MSpec}(A_{\eta})$ consists of the single point $\eta$ and hence is a dense open subscheme of $X_{MSch}$. $A_{\eta}$ is a pointed free abelian group of finite rank (hence isomorphic to some $\mathbb{Z}^n$ with base point adjoined) and we can consider the co-group morphisms 
    \begin{align}
        &m: A_{\eta} \rightarrow A_{\eta} \wedge A_{\eta}, \quad a \mapsto a \wedge a
        \\&e: A_{\eta} \rightarrow S^0=\{0, 1\}, \quad a\not = 0 \mapsto 1, 0 \mapsto 0
        \\&inv: A_{\eta} \rightarrow A_{\eta}, \quad a \mapsto a^{-1},
    \end{align}
    where $\wedge$ is the smash product (see \cite[Section 1]{monoid}). These make $\textrm{MSpec}(A_{\eta})$ a group object in the category of monoid schemes. Further, as $A$ is cancellative, $A \rightarrow A_{\eta}$ is an injection for any monoid corresponding to an affine open of $X_{MSch}$ and we can consider the co-group action
    \begin{align}
        g: A \rightarrow A \wedge A_{\eta}, \quad a \mapsto a \wedge a,
    \end{align}
    which induces a group monoid scheme action of $\textrm{MSpec}(A_{\eta})$ on $\textrm{MSpec}(A)$. This is compatible with gluing affine monoid schemes so we get an action 
    \begin{align}
    g: \textrm{MSpec}(A_{\eta}) \times X_{MSch} \rightarrow X_{MSch}.
    \end{align}
    Let $R$ be a commutative ring. As $R[S^0] = R$, $R$-realization (see \cite[Section 5]{monoid}) induces analogous co-group morphisms
    \begin{align}
        &m_R: R[A_{\eta}] \rightarrow R[A_{\eta}] \otimes_R R[A_{\eta}], \quad a \mapsto a \otimes a
        \\&e_R: R[A_{\eta}] \rightarrow R, \quad a\not = 0 \mapsto 1, 0 \mapsto 0
        \\&inv_R: R[A_{\eta}] \rightarrow R[A_{\eta}], \quad a \mapsto a^{-1},
    \end{align}
    making $\textrm{Spec} (R[A_{\eta}])$ a group object in schemes over $\textrm{Spec} (R)$. Note that $\textrm{Spec} (R[A_{\eta}])$ is a (split) torus over $R$. Then
    \begin{align}
        g_R: R[A] \rightarrow R[A] \otimes_R R[A_{\eta}], \quad a \mapsto a \otimes a
    \end{align}
    induces a group scheme action of $\textrm{Spec} (R[A_{\eta}])$ on $\textrm{Spec} (R[A])$ and $X_{MSch}$ separated implies the compatibility with gluing over monoid schemes extends to schemes. Hence we get an action on $(X_{MSch})_R$
    \begin{align}
        g_R: \textrm{Spec} (R[A_{\eta}]) \times_{\textrm{Spec} \, R} (X_{MSch})_R \rightarrow (X_{MSch})_R.
    \end{align}
    Moreover, the density of $\textrm{MSpec}(A_{\eta})$ in $X_{MSch}$ implies $\textrm{Spec} (R[A_{\eta}])$ is a dense open subscheme of $(X_{MSch})_R$. We call $T_{(X_{MSch})_R}$ an "algebraic torus over $R$".
\end{construction}

\begin{lemma}
    Let $X$ be a toric $\Lambda$-scheme. Then $X$ has a dense open group subscheme $T_X: = \textrm{Spec} \, \mathbb{Z}[A_{\eta}] = \textrm{Spec} \, \mathbb{Z}[\mathbb{Z}^n]$, an "algebraic torus over $\mathbb{Z}$", with an action
    $$
    g_{\mathbb{Z}}: T_X \times_{\mathbb{Z}} X \rightarrow X
    $$
    that restricts to the action of $T_X$ on itself. Further, the fiber of $X$ over each $(p) \in \textrm{Spec} \, \mathbb{Z}$, $X_{p} \cong (X_{MSch})_{\mathbb{F}_p}$, is a toric variety (not necessarily normal) in the sense of \cite[Section 2]{BM1} with a dense open torus given by the fiber $T_{X_{\mathbb{F}_p}} = T_{(X_{MSch})_{\mathbb{F}_p}}$ and action corresponding to the fiber of the action of $T_X$ on $X$.
\end{lemma}

\textit{Proof.} Follows by construction \ref{construction7311} when $R=\mathbb{Z},\mathbb{F}_p$. \qedsymbol


\subsection{Invariance of closed subschemes}
\label{C7_3_2} 
We describe invariant closed subschemes of toric $\Lambda$-schemes under the action of $T_X$ and its relation to fibers over $\textrm{Spec} \, \mathbb{Z}$.

\begin{observation}
    Let $X_{MSch}$ be a connected separated cancellative torsionfree monoid scheme of finite type, $R$ a commutative ring  and $i: C \hookrightarrow (X_{MSch})_R$ a closed subscheme. $C$ is invariant under the action of the torus $T_{(X_{MSch})_R}$ if the morphism 
    $$
    C \times_{R} T_{X} \xrightarrow{(i, id)} X \times_{R} T_X \xrightarrow{g_{R}} X
    $$
    factors through $C$, i.e. there exists a commuting diagram
    $$
    \begin{tikzcd}
    & C \times_{R} T_{X} \arrow[]{r}{(id, i_{R})}
    \arrow[swap]{dr}
    & X \times_{R} T_X \arrow{r}{g_{R}}
    & X
    \\
    & 
    & C \arrow{ur} 
    &  
    \end{tikzcd}
    $$
    where $C \rightarrow X$ is given by $i$. This is equivalent to showing on every affine open of $X$, $\textrm{Spec} \, R[A]$, where $T_X=\textrm{Spec} \,  R[A_{\eta}]$ and $C$ is given by an ideal $I \subset R[A]$, the corresponding diagram of rings
    $$
    \begin{tikzcd}
    & R[A]/I \otimes_{R} R[A_{\eta}] 
    & R[A] \otimes_{R} R[A_{\eta}] 
    \arrow[swap]{l}{(id, i_{R})^{\sharp}}
    & R[A]
    \arrow[swap]{l}{(g_{R})^{\sharp}}
    \arrow{dl}
    \\
    & 
    & R[A]/I \arrow{ul} 
    &  
    \end{tikzcd}
    $$
    commutes, i.e. each element of $I$ maps to $0$ under $(id, i_{R})^{\sharp} \circ (g_{R})^{\sharp}$.
\end{observation}

\begin{lemma}[Equivariant closed subschemes]
    Let $i:C \rightarrow (X_{MSch})_R$ be a closed subscheme given by the $R$-realization of an equivarant closed monoid subscheme \cite[Section 2]{monoid} then $C$ is invariant under the action of $T_{(X_{MSch})_R}$. 
\end{lemma}

\textit{Proof.} Commutativity of the affine diagram holds for $I$ generated by monoid elements as a monoid element $a \in I \mapsto a \otimes a \mapsto (a+I) \otimes a = 0 \otimes a$. \qedsymbol

\begin{example}
    Let $X$ be a toric $\Lambda$-scheme and $i:C \rightarrow X$ a closed subscheme given by the $\mathbb{Z}$-realization of an equivariant closed monoid subscheme, then $C$ is invariant under $T_X$.
\end{example}

\begin{lemma}[Invariance over $\mathbb{Z}$ implies invariance over fibers]
    Let $X$ be a toric $\Lambda$-scheme and $i_{\mathbb{Z}}:C \rightarrow X$ a $T_X$-invariant closed subscheme. Then the closed subscheme of the toric variety $X_p$, $i_{p}: C_p \rightarrow X_p$, induced by taking the fiber over a prime ideal $(p) \in \textrm{Spec} \, \mathbb{Z}$ is $T_{X_p}$-invariant.
\end{lemma}

\textit{Proof.} It suffces to restrict to the affine case, where invariance implies we have a commuting diagram as above for $R=\mathbb{Z}$. For any prime ideal $(p) \subset \mathbb{Z}$, tensoring $\otimes_{\mathbb{Z}} \mathbb{F}_p$ (with $\mathbb{F}_0=\mathbb{Q}$) we obtain a commuting diagram
$$
\begin{tikzcd}
& (\mathbb{Z}[A]/I \otimes_{\mathbb{Z}} \mathbb{Z}[A_{\eta}]) \otimes_{\mathbb{Z}} \mathbb{F}_p 
& (\mathbb{Z}[A] \otimes_{\mathbb{Z}} \mathbb{Z}[A_{\eta}]) \otimes_{\mathbb{Z}} \mathbb{F}_p  
\arrow[swap]{l}{(id, i_{\mathbb{Z}})^{\sharp} \otimes id}
& (\mathbb{Z}[A]) \otimes_{\mathbb{Z}} \mathbb{F}_p 
\arrow[swap]{l}{(g_{\mathbb{Z}})^{\sharp} \otimes id}
\arrow{dl}
\\
& 
& (\mathbb{Z}[A]/I) \otimes_{\mathbb{Z}} \mathbb{F}_p \arrow{ul}.
&  
\end{tikzcd}.
$$
By right exactness $(\mathbb{Z}[A]/I) \otimes_{\mathbb{Z}} \mathbb{F}_p \cong \mathbb{F}_p[A]/(I_{\mathbb{F}_p})$, where $I_{\mathbb{F}_p}:=im(I \otimes_{\mathbb{Z}} \mathbb{F}_p \rightarrow \mathbb{Z}[A] \otimes_{\mathbb{Z}} \mathbb{F}_p)$ is the ideal defining the fiber $C_p$. Thus we get a commuting diagram 
$$
\begin{tikzcd}
& \mathbb{F}_p[A]/(I_{\mathbb{F}_p}) \otimes_{\mathbb{F}_p} \mathbb{F}_p[A_{\eta}])  
& (\mathbb{F}_p[A] \otimes_{\mathbb{F}_p} \mathbb{F}_p[A_{\eta}])  
\arrow[swap]{l}{(id, i_{\mathbb{F}_p})^{\sharp}}
& (\mathbb{F}_p[A])  
\arrow[swap]{l}{(g_{\mathbb{F}_p})^{\sharp} }
\arrow{dl}
\\
& 
& (\mathbb{F}_p[A])/(I_{\mathbb{F}_p}) 
\arrow{ul} 
&  
\end{tikzcd},
$$
where under the isomorphisms, $(id, i_{\mathbb{Z}})^{\sharp} \otimes id$ induces $(id, i_{\mathbb{F}_p})^{\sharp}$, where $i_{\mathbb{F}_p}^{\sharp}: \mathbb{F}_p[A] \rightarrow \mathbb{F}_p[A] / (I_{\mathbb{F}_p})$ is the projection and $(g_{\mathbb{Z}})^{\sharp} \otimes id$ induces $(g_{\mathbb{F}_p})^{\sharp}$; the co-action of the algebraic torus (over $\mathbb{F}_p$). But this is precisely the description of $T_{X_p}$-invariance for $R=\mathbb{F}_p$. \qedsymbol


\subsection{Unpointed monoid morphisms}
\label{C7_3_3} 

We describe a class of monoid scheme morphisms that define the embedding in \cref{theorem7201}.

\begin{definition}
    Let $f: A \rightarrow B$ a morphism of monoids, we will say $f$ is \textit{unpointed} if $f^{-1}(0)=0$, i.e. $f$ is induced by an unpointed monoid morphism.
\end{definition}

\begin{definition}
    \label{definition7331} 
    Let $f: X_{MSch} \rightarrow M_{MSch}$ be an affine morphism of monoid schemes (see \cite[Section 6]{monoid}). $f$ is \textit{unpointed} if each induced monoid morphism $A \rightarrow B$ is unpointed. Being unpointed behaves well with respect to localizations so to show $f$ is unpointed it suffices to show it is on an affine open cover.
\end{definition}

\begin{lemma}
    \label[lemma]{lemma7333} 
    Let $f: X_{MSch} \rightarrow M_{MSch}$ be an affine morphism of connected separated cancellative torsionfree monoid scheme of finite type, e.g. a closed immersion as in \cref{theorem7201}.  Write $\eta \in X_{MSch}$, $\xi \in M_{MSch}$ for the unique minimal points and suppose $f$ was unpointed. Then for any ring $R$, taking $R$-realizations (see \cite[Section 5]{monoid}) we have a commuting diagram
    $$
    \begin{tikzcd}
	{(T_X)_R \times_R X_R} & {} & {(T_M)_R \times_R M_R} \\
	\\
	{X_R} && {M_R}
	\arrow["{g_B}", from=1-1, to=3-1]
	\arrow["{(\pi'_R, f_R)}", from=1-1, to=1-3]
	\arrow["{f_R}", from=3-1, to=3-3]
	\arrow["{g_A}", from=1-3, to=3-3]
    \end{tikzcd}.
    $$
\end{lemma}

\textit{Proof.} Let $\textrm{MSpec}(B)$ be an affine open neighbourhood of $\xi$ and $\textrm{MSpec}(A)$ an affine open neighbourhood of $\eta$ with $\textrm{MSpec}(A) = f^{-1}(\textrm{MSpec}(B))$. Then $f$ is induced by a monoid morphism $\pi: B \rightarrow A$. This morphism being unpointed implies on stalks it induces $\pi':B_{\xi} \rightarrow A_{\eta}$. Thus the following diagram (given by cogroup morphisms) commutes
$$
\begin{tikzcd}
	{B_{\xi} \wedge B} & {} & {A_{\eta} \wedge A} \\
	\\
	B && A
	\arrow["{g_B}"', from=3-1, to=1-1]
	\arrow["{(\pi', \pi)}", from=1-1, to=1-3]
	\arrow["\pi", from=3-1, to=3-3]
	\arrow["{g_A}"', from=3-3, to=1-3]
\end{tikzcd}
$$
and for a ring $R$ the $R$-realization of this diagram commutes also. This implies that the $R$-realization $f_R: (X_{MSch})_R \rightarrow (M_{MSch})_R$ restricted to the algebraic torus $(T_{X_{MSch}})_R$ of $(X_{MSch})_R$ (which is a dense open subscheme) induces a morphism $h_R:(T_{X_{MSch}})_R \rightarrow (T_{M_{MSch}})_R$ to the torus of $(M_{MSch})_R$ given by an underlying group morphism. Further, as the above diagram is compatible with gluing along monoid scheme intersections, $f_R$ will be equivariant with respect to $h_R$ and the diagram commutes. \qedsymbol

\begin{corollary}
    Let $f: X \rightarrow M$ be a closed embedding of toric $\Lambda$-schemes induced by a closed immersion $X_{MSch} \rightarrow M_{MSch}$ that is unpointed. Then $f$ is an "equivariant" embedding (as described above) and each fiber over a prime ideal $(p) \subset \mathbb{Z}$, $f_p: X_p \rightarrow M_p$ is an equivariant embedding of toric varieties over $\mathbb{F}_p$ in the sense of \cite{BM1}.
\end{corollary}

\begin{remark}
    If the unpointed monoid scheme morphism of \cref{lemma7333} is a closed embedding as in the above corollary then we further get that the morphism of group completions is surjective. It follows that for any $R$-realization the ring morphism underlying the morphism of algebraic tori is surjective. In particular, this applies to the embedding of toric $\Lambda$-schemes in \cref{theorem7201}.
\end{remark}


\section{Normal toric \texorpdfstring{$\Lambda$}{Λ}-schemes}
\label{C7_4}
We describe $\Lambda$-schemes analogous to normal toric varieties. We discuss properties associated to their fans and characterize the embedding of \cref{theorem7201} as locally given by binomial equations.

\begin{definition}
    A $\Lambda$-scheme $X$ such that $X \cong (X_{MSch})_{\mathbb{Z}}$ for $X_{MSch}$ a normal toric monoid scheme (see \cref{C7_3}) is a \textit{normal} toric $\Lambda$-scheme. In particular, a normal toric $\Lambda$-scheme has an associated fan $\Sigma$ with lattice $N$ coming from $X_{MSch}$ given by the $\mathbb{Z}$-linear dual of $A_{\eta}$ for $\eta$ the unique minimal point. Conversely, a fan $\Sigma$ will induce a normal toric $\Lambda$-scheme (see \cite[Section 4]{monoid}). A normal toric $\Lambda$-scheme $X$ will be covered by affine open integral $\Lambda$-subschemes corresponding to cones $\sigma \in \Sigma$
    $$
    U_{\sigma}:=\textrm{Spec} \, \mathbb{Z}[A(\sigma)],
    $$
    for $A(\sigma)$ a finitely generated torsionfree normal monoid.
\end{definition}

\begin{remark}
    Borger calls these $\Lambda$-schemes toric varieties (see \cite[Section 2.4]{borger2009lambdarings}). As we will be also working with usual toric varieties over fields we reserve this term for the field case.
\end{remark}

\begin{observation}
    Let $X$ be a normal toric $\Lambda$-scheme with fan $\Sigma$. Then for each prime ideal $(p) \in \textrm{Spec}  \,\mathbb{Z}$, the fiber $X_p$ is the usual toric variety over $\mathbb{F}_p$ associated to the fan $\Sigma$ (see \cite[Chapter 3]{cox2011toric}). Further, the algebraic torus and action described in the previous section given by the fiber of $T_X$ on $X$ is precisely the usual torus and action given by the character lattice (see \cite[Chapter 1, Chapter 3]{cox2011toric}, \cite[Section 2]{BM1}).
\end{observation}


\subsection{Smooth normal toric \texorpdfstring{$\Lambda$}{Λ}-schemes}
\label{C7_4_1} 
Let $X$ be a normal toric $\Lambda$-scheme associated to a fan $\Sigma$ and normal toric monoid scheme $X_{MSch}$. We will say $X$ is smooth if it is smooth over $\textrm{Spec} \, \mathbb{Z}$, which is equivalent to $X_{MSch}$ being a smooth monoid scheme (and $\Sigma$ being regular) by \cref{lemma103}.  
\\
\\Let $X$ be a smooth normal toric $\Lambda$-scheme. We state some facts analogous to those for normal toric varieties over fields. The affine opens corresponding to cones $\sigma = Cone(e_1, \dots, e_r) \in \Sigma$ are of the form
$$
U_{\sigma}:=\textrm{Spec} \, \mathbb{Z}[A(\sigma)] = \textrm{Spec} \, \mathbb{Z}[x_1, \dots, x_r, y_{r+1}^{\pm}, \dots, y_n^{\pm}],
$$
where $n$ is the dimension of the lattice. The torus in this affine open is
$$
T_X = \textrm{Spec} \, \mathbb{Z}[x_1^{\pm}, \dots, x_r^{\pm}, y_{r+1}^{\pm}, \dots, y_n^{\pm}].
$$
Let $\Delta \leq \sigma$ be a face. Then $\Delta=Cone(e_{j_1}, \dots, e_{j_s})$ for some subset $\{ e_{j_1}, \dots, e_{j_s}\} \subseteq \{e_1, \dots, e_r\}$ of the vertices of $\sigma$. Write $V_{\Delta} :=\{j_1, \dots, j_s  \} \subseteq \{1, \dots, r\}$. We have an equivariant closed subscheme $X_{MSch}(Star(\Delta))$ of $X_{MSch}$ given by the fan $Star(\Delta)$ (see \cite[Exercise 3.2.7]{cox2011toric}). Then we get a closed subscheme of $X$
$$
Z_{\Delta} := (X_{MSch}(Star(\Delta)))_{\mathbb{Z}},
$$
which on each affine open $U_{\sigma}$ with $\Delta \leq \sigma \in \Sigma$ is given by the ideal generated by the monoid elements $x_i, i \in V_{\Delta}$. In particular, $Z_{\Delta}$ will be $T_X$-invariant and a closed $\Lambda$-subscheme (with toric $\Lambda$-structure) as it is generated by monoid elements. Further, $Star(\Delta)$ is a regular fan so $Z_{\Delta}$ itself is a smooth normal toric $\Lambda$-scheme. The fibers $(Z_{\Delta})_p$ will precisely correspond to the usual smooth torus-invariant subschemes given by the orbit closures of the toric varieties $X_p$.

\begin{example}[normal crossings divisors]
    \label[example]{example7412} 
    Let $X$ be a smooth normal toric $\Lambda$-scheme with fan $\Sigma$. Let $E=\sum_i E^i$, where $i$ sums across the (finitely many) cones $\Delta \in \Sigma$ with $E^i = Z_{\Delta}$ of codimension one in $X$. As  the $x_i, i=1, \dots, r$ form a regular sequence in $\mathbb{Z}[x_1, \dots, x_r, y_{r+1}^{\pm}, \dots, y_n^{\pm}]$ they will remain a regular sequence in the stalk at any $x \in U_{\sigma}$. It follows that as $Z_{\Delta}$ is defined by a $x_i$ in $U_{\sigma}$ that $E$ is a $\mathbb{Z}$-flat simple normal crossings divisor (see \cite[Section 2.3]{OrderReductionLambdamarkedMonomialIdeals}). Moreover, $E$ is locally toric (see \cite[Definition 2.3.14]{OrderReductionLambdamarkedMonomialIdeals}) for the affine $\Lambda$-cover of $X$ given by the cones of $\Sigma$, which is induced from the corresponding affine monoid scheme cover.
\end{example}

\begin{remark}
    \label[remark]{remark7413} 
     $Z_{\Delta}$ will meet the conditions of a blow-up centre that induces a $\Lambda$-structure as in \cite[Section 1]{OrderReductionLambdamarkedMonomialIdeals}. We can see this as $Z_{\Delta}$ is the $\mathbb{Z}$-realization of an equivariant closed monoid subscheme of $X_{MSch}$. Alternatively, $Z_{\Delta}$ is a finite intersection $\cap_J E^j$ for the locally toric simple normal crossings divisor $E$ above so we can use 
     \cite[Example 2.3.5]{OrderReductionLambdamarkedMonomialIdeals}. This second perspective also shows that $Z_{\Delta}$ has normal crossings with $E$ and so we can take total transforms as in \cite[Definition 2.3.7]{OrderReductionLambdamarkedMonomialIdeals}.
\end{remark}


\subsection{Binomial description of the closed embedding}
\label{C7_5} 
We now describe the embedding of \cref{theorem7201} locally as quotients of polynomial rings over $\mathbb{Z}$ by binomial equations. Let $X = \textrm{Spec}(\mathbb{Z}[A])$ be an affine toric $\Lambda$-scheme with $A$ a cancellative torsionfree monoid of finite type. It follows from the same arguments for toric varieties over fields (see \cite[Proposition 1.1.9]{cox2011toric}) that

$$
\mathbb{Z}[A] \cong \mathbb{Z}[x_1, \dots, x_s]/ I,
$$
where $I$ is prime and generated by binomials.
\\
\\Let $B \rightarrow A$ be a surjective monoid morphism of torsionfree cancellative monoids of finite type, inducing a surjective morphism of monoid algebras $\mathbb{Z}[B] \rightarrow \mathbb{Z}[A]$. Suppose $B=\langle x_1, \dots, x_r, y_{r+1}^\pm, \dots, y_n^{\pm}\rangle$, i.e. the monoid corresponding to a regular cone. Then 
\begin{align}
  \mathbb{Z}[B] &\cong \mathbb{Z}[x_1, \dots, x_r, {y_{r+1}}^{\pm }, \dots, {y_n}^{\pm }]
  \label{eqn7.11}
  \\&\cong \frac{\mathbb{Z}[x_1, \dots, x_r, {y_{r+1}}, \dots, {y_n}, u_{r+1}, \dots, u_{n}]}{(u_iy_i-1, i=r+1, \dots ,n)} .
  \label{eqn7.12} 
\end{align}
Now by surjectivity we can take the image of the generators of $B$ to be the generators of $A$. In particular, there are $2n-r$ of them and $2(n-r)$ (the $y,u$-variables) are units and hence their images in $A$ are also. Thus, by the previous arguments on monoid algebras given by binomials we can write 
\begin{align}
  \mathbb{Z}[A] \cong \frac{\mathbb{Z}[x_1, \dots, x_r, {y_{r+1}}, \dots, {y_n}, u_{r+1}, \dots, u_{n}]}{I}  
  \label{eqn7.13}
\end{align}
for $I$ generated by binomials, where some of these binomials are $u_iy_i-1, i=r+1, \dots ,n$ making $y_i, u_i$ units. In particular, we can write this as
\begin{align}
   \mathbb{Z}[A] \cong \frac{\mathbb{Z}[x_1, \dots, x_r, {y_{r+1}}^{\pm}, \dots, {y_n}^{\pm}]}{I}, 
   \label{eqn7.14} 
\end{align}
where the binomials generating $I$ may feature laurent monomials in the $y$-variables.

\begin{remark}
    Recall a surjective monoid morphism is given by the quotient by a congruence. We can interpret these binomials defining $\mathbb{Z}[A]$ in $\mathbb{Z}[B]$ as the quotient by the "$\mathbb{Z}$-realization" of this congruence. 
\end{remark}

\begin{observation}
    \label{observation7522} 
    Let $X$ be a toric $\Lambda$-scheme and $M$ a smooth toric $\Lambda$-scheme with fan $\Sigma$. Suppose $X \hookrightarrow M$ is a closed embedding induced by a closed immersion $X_{MSch} \rightarrow M_{MSch}$ of the underlying monoid schemes. For each cone $\sigma \in \Sigma$ and affine open subscheme $U_{\sigma}$ given by a ring of the form (\ref{eqn7.11}, \ref{eqn7.12}), $X \cap U_{\sigma}$ is given by rings of the form (\ref{eqn7.13}, \ref{eqn7.14}) and as a closed embedding is given by an ideal generated by binomials. Further, each embedding of toric varieties given by fibers over $(p) \in \textrm{Spec}  \, \mathbb{Z}$, $X_p \hookrightarrow M_p$, sees $X_p$ locally cut out by precisely the same binomials (see \cite[Section 2]{BM1}).
\end{observation}

\textbf{The torus.} Let $M$ be a smooth toric $\Lambda$-scheme and $U_{\sigma}$ an affine open chart corresponding to the monoid $B=\langle x_1, \dots, x_r, y_{r+1}^\pm, \dots, y_n^{\pm}\rangle$. Then the torus $T_M$ is given by the $\mathbb{Z}$-realization of the localization $B_{x_1 \cdots x_r} = \langle x_1^\pm, \dots, x_r^{\pm}, y_{r+1}^\pm, \dots, y_n^{\pm}\rangle$, i.e. in binomials we have
$$
T_M = \textrm{Spec} \, \mathbb{Z}[x_1^\pm, \dots, x_r^{\pm}, y_{r+1}^\pm, \dots, y_n^{\pm}] \subset \textrm{Spec} \,  \mathbb{Z}[x_1, \dots, x_r, y_{r+1}^\pm, \dots, y_n^{\pm}].
$$
If we further assume $X \hookrightarrow M$ comes from an unpointed embedding (see \cref{C7_3_3}) then locally the images of the generators of $B$ are not $0$ in the monoid $A$ underlying $U_{\sigma} \cap X$. In particular, this means the torus $T_X$, given by inverting the generators of $A$, which are the images of the generators of $B$ (which may be sent to $1$ or identified with each other but never $0$) has the form
$$
T_X = \textrm{Spec} \, \frac{\mathbb{Z}[x_1^\pm, \dots, x_r^{\pm}, y_{r+1}^\pm, \dots, y_n^{\pm}]}{I} \subset \textrm{Spec} \, \frac{\mathbb{Z}[x_1, \dots, x_r, {y_{r+1}}^{\pm}, \dots, {y_n}^{\pm}]}{I}.
$$

\textbf{Orbit closures in binomials.} Recall from \cref{C7_4_1} the equivariant closed subschemes $Z_{\Delta}$ (orbit closures) of smooth toric $\Lambda$-schemes $M$, where $\Delta \in \Sigma$ is a cone of the fan. Then for $\sigma=Cone(e_1, \dots, e_r) \in \Sigma$ with $\Delta$ a face, $\Delta$ is given by its vertex set 
$$
V_{\Delta}, \quad \Delta = Cone(e_i:i \in V_{\Delta}) \subset \sigma.
$$
Then $Z_{\Delta}$ as a closed subscheme of $M$ will be locally cut out in a $U_{\sigma}$ by $I=(x_i: i \in V_{\Delta})$, i.e. 
$$
Z_{\Delta} \cap U_{\sigma} \cong \textrm{Spec} \, \frac{\mathbb{Z}[x_1, \dots, x_r, {y_{r+1}}^{\pm}, \dots, {y_n}^{\pm}]}{(x_i: i \in V_{\Delta})}.
$$
\\\textbf{Smooth embeddings of toric $\Lambda$-schemes.} Let $X \hookrightarrow M$ be as in \cref{theorem7201}. If $X$ itself is a smooth normal toric $\Lambda$-scheme corresponding to a fan $\Sigma'$ then $U_{\sigma} \cap X$ will be an affine open $U_{\tau}$ of $X$ corresponding to a regular cone $\tau \in \Sigma'$. It follows that
\begin{align*}
    U_{\tau} &= \textrm{Spec} \, \mathbb{Z}[x_1, \dots, x_s, {y_{s+1}}^{\pm }, \dots, {y_m}^{\pm }] 
    \\&\cong X \cap U_{\sigma} = \textrm{Spec} \, \frac{\mathbb{Z}[x_1, \dots, x_r, {y_{r+1}}^{\pm}, \dots, {y_n}^{\pm}]}{I}, \quad s \leq r.
\end{align*}

Now $X \hookrightarrow M$ coming from an unpointed closed embedding implies these morphisms are induced by a surjective morphism of the lattices with $\tau$ mapping into $\sigma$. We can then choose a basis of the lattice of $X$ extending the generators of $\tau$ and extend to a basis of the lattice of $M$ given by $\sigma$. Thus we can assume that $x_1, \dots, x_s$ come from the image of a subset of the $x_1, \dots, x_r$.


\section{Blow-ups and strict transforms}
\label{C7_6}
Let $X$ be a toric $\Lambda$-scheme and $M$ a smooth toric $\Lambda$-scheme with fan $\Sigma$. Let $X \hookrightarrow M$ be a closed embedding induced by a closed immersion $X_{MSch} \rightarrow M_{MSch}$ of the underlying monoid schemes. We describe the effect of blowing up $M$ along a $Z_{\Delta}$ corresponding to $\Delta \in \Sigma$. We consider the strict transform of $X \hookrightarrow M$ and global properties induced by the underlying monoid schemes.


\subsection{Blowing up}
\label{C7_6_1} 
Let $\sigma \in \Sigma$ be a maximal cone with $\Delta$ as a face and write $I = (x_i: i \in V_{\Delta}) \subset \mathbb{Z}[x_1, \dots, x_r, {y_{r+1}}^{\pm}, \dots, {y_n}^{\pm}]$ for the ideal defining $Z_{\Delta} \cap U_{\sigma}$. For each $i \in V_{\Delta}$ the blow-up $B_{Z_{\Delta}} M$ has an affine open $U_{\sigma_i}$ given by the affine blow-up algebra corresponding to $x_i$. Then
$$
U_{\sigma_i} \cong \textrm{Spec} \,  \mathbb{Z}[w_1, \dots, w_r, {y_{r+1}}^{\pm}, \dots, {y_n}^{\pm}]
$$
and $U_{\sigma_i} \rightarrow U_{\sigma}$ is induced by 
$$
\mathbb{Z}[x_1, \dots, x_r, {y_{r+1}}^{\pm}, \dots, {y_n}^{\pm}] \rightarrow \mathbb{Z}[w_1, \dots, w_r, {y_{r+1}}^{\pm}, \dots, {y_n}^{\pm}]
$$
given by the standard substitution rule
\begin{align}
    &y_k \mapsto y_k, \quad k=r+1, \dots, n
    \\&x_j \mapsto w_iw_j, \quad j \in V_{\Delta} \backslash \{i\}
    \\&x_j \mapsto w_j, \quad j \not \in V_{\Delta} \backslash \{i\}.
\end{align}
Moreover, these $U_{\sigma_i}$ cover $B_{Z_{\Delta}} M$ as $\sigma$ varies along maximal cones of $\Sigma$.

\begin{remark}
    $B_{Z_{\Delta}} M$ is a smooth normal toric $\Lambda$-scheme corresponding the fan $\Sigma^{\ast}(\Delta)$ given by star subdivision. Then each $U_{\sigma_i}$ corresponds to a maximal cone $\sigma_i = Cone(e_0, e_1, \dots, \hat{e_i}, \dots, e_r) \in \Sigma^{\ast}(\Delta)$, where $e_0$ is the barycentre (see \cite[Section 11]{monoid}).
\end{remark}



\subsection{Strict transforms}
\label{C7_6_2}  
Let $\sigma \in \Sigma$ be a maximal cone with $\Delta$ as a face and let $x^{\alpha}-x^{\beta}y^{\gamma}$ be a binomial generating $I \subset \mathbb{Z}[x_1, \dots, x_r, {y_{r+1}}^{\pm}, \dots, {y_n}^{\pm}]$, where $I$ defines $X \cap U_{\sigma}$ (see \cref{C7_5}). Note we can write binomials in this form by moving all $y$-variables to one monomial via multiplication by inverses. Write $\alpha=(\alpha_1, \dots, \alpha_r), \beta=(\beta_1, \dots, \beta_r) \in \mathbb{N}^r$ and $\gamma = (\gamma_1, \dots, \gamma_{n-r}) \in \mathbb{Z}^{n-r}$. The total transform of $X$ in the affine open $U_{\sigma_i}$ will be generated by applying the substitution rule to the binomials of $I$, i.e. $x^{\alpha}-x^{\beta}y^{\gamma}$ becomes
$$
w^{\alpha'}-w^{\beta'}y^{\gamma}, \quad \alpha'=(\alpha_1', \dots, \alpha_r') \quad  \beta'=(\beta_1', \dots, \beta_r'),
$$
where $\alpha_j' = \alpha_j, \beta_j'=\beta_j$ for all $j \not = i$ and 
$$
\alpha_i' = \alpha_\Delta: = \sum_{j \in V_{\Delta}} \alpha_j, \quad \beta_i' = \beta_\Delta: =\sum_{j \in V_{\Delta}} \beta_j.
$$
Then writing $I'$ for the ideal generated by the binomials $x^{\alpha'}-x^{\beta'}y^{\gamma}$, $\tilde{X} \cap U_{\sigma_i}$ is   
$$
\textrm{Spec} \left(\frac{\mathbb{Z}[w_1, \dots, w_r, {y_{r+1}}^{\pm}, \dots, {y_n}^{\pm}]}{I'} \right)
$$
for $\tilde{X}$ the total transform of $X$. Writing $R$ for the underlying ring of $\tilde{X} \cap U_{\sigma_i}$ we can compute the strict transform $X'$ from the total transform. $X' \cap U_{
\sigma_i}$ is given by the quotient $\frac{R}{J}$, where 
$$
J: = ker\{ R \rightarrow R_{w_i}    \} = \{  r \in R: w_i^nr = 0 \textrm{ for some }n > 0       \}.
$$

\begin{example}[Strict transform of affine hypersurfaces]
    Let $X,M$ be as above and both affine and let $X \hookrightarrow M$ be a hypersurface. Then $\Sigma$ corresponds to a single regular cone and its faces and $I$ as above is generated by a single binomial $x^{\alpha}-x^{\beta}y^{\gamma}$. Assume that $\alpha_\Delta \not = 0, \beta_\Delta \not = 0$ and without loss of generality assume $\alpha_\Delta \leq \beta_\Delta$ and write $d_{\Delta} = \beta_\Delta-\alpha_\Delta$. Then $I'$ corresponding to $\tilde{X} \cap U_{\sigma_i}$ is generated by
    $$
    w^{\alpha'}-w^{\beta'}y^{\gamma} = w_i^{\alpha_{\Delta}}(w^{\alpha''}-w^{\beta''}y^{\gamma}), 
    $$
    where $\alpha_j'' = \alpha_j', \beta''_j= \beta_j'$ for $ j \not = i$ and  $\alpha_i'' = 0, \beta_i'' = d_{\Delta}$. Then
    $$
    X' \cap U_{\sigma_i} \cong \textrm{Spec} \left(    \frac{\mathbb{Z}[w_1, \dots, w_r, {y_{r+1}}^{\pm}, \dots, {y_n}^{\pm}]}{(w^{\alpha''}-w^{\beta''}y^{\gamma})}          \right).
    $$
\end{example}

\begin{remark}
    In \cref{C7_8} on the resolution of a hypersurface we will blow up with centres $Z_{\Delta}$ such that both $\alpha_{\Delta}, \beta_{\Delta}$ are nonzero and so these local descriptions of the strict transforms apply. We can also consider the situation where one of $\alpha_{\Delta}, \beta_{\Delta}$ is 0. In this case we can use similar arguments to show $J=0$ and the local strict transform is the local total transform.
\end{remark}



\subsection{Blow-up properties of monoid schemes}
\label{C7_7}  
Let $X \hookrightarrow M$ be a closed embedding of toric $\Lambda$-schemes as in \cref{theorem7201}, induced by $X_{MSch} \rightarrow M_{MSch}$ with $\Sigma$ the regular fan corresponding to $M$. Let $X'$ be the strict transform of $X$ given by the blow-up of $M$ with $Z_{\Delta}$ and $M'$ the blow-up itself. We have seen from \cref{C7_4_1} that $Z_{\Delta}$ given by $\Delta \in \Sigma$ is itself a smooth normal toric $\Lambda$-scheme and as a closed subscheme is $\Lambda$-equivariant and $T_M$-invariant. Thus, condition 1 of \cref{theorem7201} is always satisfied for the blow-ups we consider. Further, $M'$ corresponds to a regular fan given by star-subdivision, i.e. the ambient scheme in any of the steps in a blow-up sequence is always a smooth toric $\Lambda$-scheme. The strict transform $X' \hookrightarrow M'$ also retains the same structure as the initial embedding of \cref{theorem7201}. $X'$ is a toric $\Lambda$-scheme and $X' \hookrightarrow M'$ is a $\Lambda$-equivariant closed embedding induced by an unpointed closed monoid subscheme. These conditions solely concern monoid schemes and the arguments follow from local monoid arguments as in \cite[Proposition 9.1, Proposition 9.2]{monoid}. Thus, the blow-ups, strict transforms and embeddings using centres of the form $Z_{\Delta}$ always retain the same structure as in the hypotheses of \cref{theorem7201}. This allows us to make use of inductive arguments.


\subsection{Smooth invariant centres}
\label{C7_6_4} 
We show that centres satisfying condition 1 of \cref{theorem7201} will be either a $Z_{\Delta}$ or a disjoint union of $Z_{\Delta}$'s and this holds over each fiber. This will apply to all blow-ups we consider. 

\begin{lemma} 
    \label[lemma]{lemma7731}
    Let $I \subset \mathbb{Z}[x_1, \dots, x_r, y_1^{\pm}, \dots, y_{n-r}^{\pm}] =R$ be an ideal such that $\frac{R}{I}$ is flat over $\mathbb{Z}$. Denote by $\mathbb{Q}I \subset\mathbb{Q}[x_1, \dots, x_r, y_1^{\pm}, \dots, y_{n-r}^{\pm}]$ the ideal generated by $I$. Then $\mathbb{Q}I \cap R = I$.
\end{lemma}

\textit{Proof.} One direction is clear. We show $\mathbb{Q}I \cap R \subset I$. Let $f \in \mathbb{Q}I \cap R$ and suppose $f \not \in I$, then $f \in \mathbb{Q}I$ implies there exists some positive integer $n$ such that $nf \in I$. Take $n_0$ to be the greatest common divisor of all such $n$. Consider $f+I$ in $R/I$, which is $n_0$-torsion. Now $R$ is torsionfree over $\mathbb{Z}$ and by $D$ being flat over $\mathbb{Z}$, $R/I$ is also. But $nf+I$ is $0+I$ so $f \in I$. \qedsymbol.

\begin{proposition}
    Let $D$ be a closed subscheme of a smooth toric $\Lambda$-scheme $M$ with fan $\Sigma$. If $D$ is smooth over $\textrm{Spec} \, \mathbb{Z}$ and $T_M$-invariant then 
    \begin{enumerate}
        \item Each fiber over $(p)$, $D_p$, is a disjoint union of orbit closures of $M_p$ and the cones $\Delta \in \Sigma$ corresponding to the orbit closures are independent of $p$.
        \item $D$ itself is a disjoint union of $Z_{\Delta}$'s of $M$ given by the same cones.
    \end{enumerate}
\end{proposition}

\textit{Proof.} By the results on $T_M$-invariance and smoothness over $\textrm{Spec} \, \mathbb{Z}$ each $D_p$ is a closed smooth invariant subvariety of $M_p$ (a smooth toric $\mathbb{F}_p$-variety) and so it follows from facts on normal toric varieties (see \cite[Section 4]{BM1}, \cite[Chapter 3]{cox2011toric}) that $D_p$ is a disjoint union of orbit closures of $M_p$, which correspond to cones. In particular, if $(U_{\sigma})_p$ is an affine open of $M_p$ corresponding to a cone of $\Sigma_{M_p}=\Sigma$ then assuming non-empty intersection $D_p \cap (U_{\sigma})_p$ is defined by the vanishing of some non-invertible variables (see \cref{C7_3_2}). Thus to prove (1) it remains to show that the cones are all the same. It suffices to show that on any affine open chart $U_{\sigma}$ of $M$ intersecting $D$ the fibers are given by the vanishing of the same variables (which will define the orbit closures) and that these variables generate the ideal defining $D$.
\\Let $U_{\sigma}$ be an affine open subscheme of $M$ corresponding to $\sigma \in \Sigma$ and consider $D \cap U_{\sigma}$. Now $D$ smooth over $\textrm{Spec} \, \mathbb{Z}$ implies $D \cap U_{\sigma}$ is smooth over $\textrm{Spec} \, \mathbb{Z}$, i.e. smooth of some relative dimension in each neighbourhood. But $(D \cap U_{\sigma})_p = D_p \cap (U_{\sigma})_p$ is an affine open of an orbit closure of $(M_{\sigma})_p$ and so is an irreducible $\mathbb{F}_p$-variety. Thus its dimension is given by transcendence degree and can be computed on any open. Therefore, the dimensions of $(D \cap U_{\sigma})_p $ are equal for all $p$.
\\Let $I \subset \mathbb{Z}[x_1, \dots, x_r, y_1^{\pm}, \dots, y_{n-r}^{\pm}] =R$ be the ideal defining $D \cap U_{\sigma}$ so $\frac{R}{I}$ is flat over $\mathbb{Z}$ by the flatness of $D$. $D_0$ (the generic fiber) is given in $(U_{\sigma})_0$ by the ideal $\mathbb{Q}I \subset\mathbb{Q}[x_1, \dots, x_r, y_1^{\pm}, \dots, y_{n-r}^{\pm}]$ and as $\mathbb{Q}I$ defines a smooth orbit closure it is generated by some $x_i$. Thus by \cref{lemma7731}, each $x_i \in I$. Now for each $x_i$ generating $\mathbb{Q}I$, $x_i \in I$ implies $x_i \in I/pI$, the ideal defining $D_p$ in $(U_{\sigma})_p$. But these $I/pI$ also define orbit closures and so are generated by $x_j$, the number of which determine the dimension. If there were any $x_j \not = x_i$ then the dimension of $D_p \cap (U_{\sigma})_p$ would be less than that of $D_0 \cap (U_{\sigma})_0$, contradicting the fact that the fibers have the same dimension. This proves (1). For (2) write $N=(x_1, \dots, x_j) \subset I=M \subset R$ for the ideal generated by the variables. Then $N \otimes_{\mathbb{Z}} \mathbb{Q} \cong M \otimes_{\mathbb{Z}} \mathbb{Q}$ and $\frac{N}{pN} \cong \frac{M}{pM}$. We show the $R$-module $\frac{M}{N}$ is $0$, which will imply $I$ defining $D$ in $U_{\sigma}$ is given by the vanishing of the same variables as the orbit closures defining $D_p$ in $(U_{\sigma})_p$. Thus $D$ is covered by affine opens, which in $M$, match the cone and fan structure of the orbit closures in each fiber and so $D$ is given by a disjoint union of "orbit closures" ($Z_{\Delta}$'s) in $M$. By hypothesis $\frac{M}{N} \otimes_{\mathbb{Z}} \mathbb{Q} \cong \frac{M \otimes_{\mathbb{Z}} \mathbb{Q}}{N \otimes_{\mathbb{Z}} \mathbb{Q}} = 0$, which means that $\frac{M}{N}$ has $\mathbb{Z}$-torsion, i.e. every element is annihilated by some integer. But $\frac{M}{N}$ is finitely generated as an $R$-module, say by $f_1, \dots, f_t$, so we can find integers $n_i$, which annihilate $f_i$. Write $n_0$ for the product of these $n_i$ so that $n_0$ annihilates every $f_i$ and so every element of $\frac{M}{N}$. Then $\frac{\frac{M}{N}}{n_0\frac{M}{N}} \cong \frac{M}{N}$. But for all $p$, $\frac{\frac{M}{N}}{p\frac{M}{N}} = 0$, which implies for all $n \geq 1$, $n\frac{M}{N}=\frac{M}{N}$. In particular, taking $n=n_0$ implies $\frac{M}{N} = 0$. \qedsymbol


\section{Resolution of a hypersurface}
\label{C7_8} 
Let $X \hookrightarrow M$ be a closed embedding of a toric $\Lambda$-scheme as in \cref{theorem7201}. In this section we prove \cref{theorem7201} in the special case that $X$ is given by a locally principal ideal, locally generated by a single binomial. The arguments follow \cite[Section 5]{BM1} though instead of using the locus of maximal order we show we can blow up components of the closed subscheme given by Hasse derivatives (which has a more obvious $\Lambda$-scheme theoretic description). We do this by showing that using Hasse derivatives give analogous facts to \cite[Section 5]{BM1}.

\subsection{Affine case}
\label{C7_8_1} 
We consider the affine case where $M$ corresponds to a single regular cone and its faces and $X$ is generated by a single binomial $f=x^{\alpha}-x^{\beta}y^{\gamma}$. Then
$$
M = \textrm{Spec} \, \mathbb{Z}[x_1, \dots, x_r, {y_{r+1}}^{\pm}, \dots, {y_n}^{\pm}], \quad X = \textrm{Spec} \left(\frac{\mathbb{Z}[x_1, \dots, x_r, {y_{r+1}}^{\pm}, \dots, {y_n}^{\pm}]}{x^{\alpha}-x^{\beta}y^{\gamma}} \right) 
$$
and
$$
 Z_{\Delta} = \textrm{Spec} \left(\frac{\mathbb{Z}[x_1, \dots, x_r, {y_{r+1}}^{\pm}, \dots, {y_n}^{\pm}]}{(x_i: i \in V_{\Delta})} \right).
$$
Without loss of generality assume $d:=|\alpha| \leq |\beta|$. Further, as $f$ is a generating binomial we can assume that $x^{\alpha}, x^{\beta}$ share no variables so we may write $x_1, \dots, x_r=u_1, \dots, u_k, v_1, \dots, v_{r-k}$ and 
$$
x^{\alpha}-x^{\beta}y^{\gamma} = u^{\alpha}-v^{\beta}y^{\gamma} \in \mathbb{Z}[u_1, \dots, u_k, v_1, \dots, v_{r-k}, {y_{r+1}}^{\pm}, \dots, {y_n}^{\pm}]
$$ 
with $\alpha \in \mathbb{N}^k, \beta \in \mathbb{N}^{r-k}$. 


\subsubsection{Binomial Hasse derivatives and the Hasse derivative subscheme}
\label{C7_8_1_2} 
We define the Hasse derivatives of a binomial as in \cite[Remark 3.7]{BM1} and describe closed subschemes defined by their underlying monomials.
\\
\\Let $f=u^{\alpha}-v^{\beta}y^{\gamma}$ as above and write $x=(u, v) = (u_1, \dots, u_k, v_1, \dots, v_{r-k}), X=(X_1, \dots, X_r) = U_1, \dots, U_k, V_1, \dots, V_{r-k}$ and $ y=(y_{r+1}, \dots, y_n)$ then 
\begin{align*}
    f(x+X, y^{\pm}) &= f(u+U,v+V, y^{\pm}) 
    \\&= (u_1+U_1)^{\alpha_i} \cdots (u_k+U_k)^{\alpha_k}-(v_1+V_1)^{\beta_1} \cdots (v_{r_k}+V_{r-k})^{\beta_{r-k}} y^{\gamma}
    \\&=\sum_{\zeta \leq \alpha}c_{\zeta}(u)^{\alpha-\zeta}(U)^{\zeta}-y^{\gamma}\sum_{\delta \leq \beta}c_{\delta}(v)^{\beta-\delta}(V)^{\delta}
\end{align*}
by multi-binomial theorem, where by lexicographic ordering $\zeta \leq \alpha$ means $\zeta_i \leq \alpha_i, i = 1, \dots, k$. For each $\zeta$ (resp. $\delta$), $c_{\zeta}(u)^{\alpha-\zeta}$ (resp. $c_{\delta}(v)^{\beta-\delta}$)  is the \textit{Hasse derivative of order $\zeta$} (resp. $\delta$) (cf. \cite[Remark 3.7]{BM1}). Note the coefficients will be the product of binomial coefficients and the Hasse derivative of order $\bar{0}$ will be $f$ itself. We will call the $u^{\alpha-\zeta}=u_1^{\alpha_1-\zeta_1} \cdots u_k^{\alpha_k-\zeta_k}$ (resp. $v^{\beta-\delta}$) the \textit{Hasse derivative monomial of order $\zeta$ (resp. $\delta$)}. Equivalently, we can consider the Hasse derivatives of $\mathbb{Z}[x_1, \dots, x_r, {y_{r+1}}^{\pm}, \dots, {y_n}^{\pm}]$ over $\mathbb{Z}$ with respect to $x_1, \dots, x_r$ in the sense of \cite{DifferentialOperatorsAlgorithmicWeightedResolution} (cf. \cite[3.74.6]{10.2307/j.ctt7rptq}). Viewing $\zeta, \delta$ as $(\zeta, \bar{0}), (\bar{0}, \delta) \in \mathbb{N}^r$ the Hasse derivative of $f$ of order $\zeta$ is $D^{\zeta}(f)$ by \cite[Lemma 7]{DifferentialOperatorsAlgorithmicWeightedResolution} with coefficient $\binom{\alpha}{\zeta}$. The following lemma compares Hasse derivatives in fibers to the pullback of Hasse derivative monomials.

\begin{lemma}
    \label[lemma]{lemma7711} 
    Let $J, I \subset \mathbb{Z}[x_1, \dots, x_r, {y_{r+1}}^{\pm}, \dots, {y_n}^{\pm}]$ be the ideals generated by all Hasse derivatives (resp. Hasse derivative monomials) of $f$ of order $< d$ (and $f$ itself) and write 
    $$
    I_p, J_p \subset \mathbb{F}_p[x_1, \dots, x_r, {y_{r+1}}^{\pm}, \dots, {y_n}^{\pm}]
    $$
    for the pullback of $I, J$ as ideals to the fiber over $(p)$. We can take the Hasse derivatives of $f$ in $\mathbb{F}_p[x_1, \dots, x_r, {y_{r+1}}^{\pm}, \dots, {y_n}^{\pm}]$ using the same definition and $J_p$ will be the ideal generated by all Hasse derivatives (with coefficient) of order $<d$. Then 
    $$
    V(I_p)=V(J_p) \subset \textrm{Spec} \, \mathbb{F}_p[x_1, \dots, x_r, {y_{r+1}}^{\pm}, \dots, {y_n}^{\pm}].
    $$
\end{lemma}
 
\textit{Proof.} For each Hasse derivative with coefficient not divisible by $p$ we can divide off the coefficient. Then $J_p$ is generated by a subset of Hasse derivative monomials of order $< d$ (and $f$) and $J_p \subseteq I_p$. Now, let $\mathfrak{p} \in V(J)$ and $v^{\beta-\delta}$ be a Hasse derivative monomial of order $< d$. Then $|\delta| < d$ and we can write $\delta=\delta'+\delta''$ where $\delta'=\delta_{1}', \dots, \delta_{r-k}'$ ($\delta'$ can be $0$) and 
$$
\delta_i' = 
\begin{cases}
    \beta_i & \textrm{if }  \delta_i = \beta_i\\
    0 & otherwise
\end{cases}
\implies
\beta_i-\delta_i' = 
\begin{cases}
    0 & \textrm{if } \beta_i-\delta_i = 0\\
    \beta_i & \textrm{if } \beta_i-\delta_i \not = 0
\end{cases}
$$
Then $|\delta'| < d$ and $v^{\beta-\delta'}$ is a Hasse derivative of order $< d$ (its coefficient is $1$) and so $v^{\beta-\delta'} \in \mathfrak{p}$, which implies $v_i \in \mathfrak{p}$ for some $i$ with $\beta_i-\delta_i \not = 0$ and so $v^{\beta-\delta} \in \mathfrak{p}$. The $u$-monomial case can be argued in a similar way. \qedsymbol


\begin{definition}
    Define $S_f(0)$ as the closed subscheme of $M$ given by the ideal generated by $f$ and Hasse derivative monomials of order $< d$. Then $S_f(0)$ is a closed subscheme of $X$ where all $u$-variables $u_i$ with $\alpha_i > 0$ vanish. 
\end{definition}

\begin{remark}
    We can view $S_f(0)$ as the $\mathbb{Z}$-realization of a closed monoid subscheme of both $X_{MSch}$ and $M_{MSch}$, equivariant closed in the case of $X_{MSch}$. Further, the preceding lemma shows the closed subscheme, $S_f(0)_p$ in $M_p$, given by a fiber over $(p)$ defines the same closed subspace in $M_p$ as the closed subscheme defined by Hasse derivatives.
\end{remark}

We state the same combinatorial conditions on $\Delta$ as in \cite[Section 5]{BM1} used for the resolution. The arguments are straightforward and follow from simple derivative properties and the arguments of \cite[Section 5]{BM1} so we omit their details. We assume $d>1$ as if $d \leq 1$ then $X$ will already be smooth (the $d=1$ case is clear and we argue the $d=0$ case later on).
\\
\\Let $Z_{\Delta}$ be as above. $Z_{\Delta} \subset S_f(0)$ if and only if 
\begin{align}
    &i \in V_\Delta \textrm{ if } \alpha_i > 0 \quad (\textrm{i.e. } \alpha_{\Delta} = d),
    \label{eqn718}
    \\&\beta_{\Delta} \geq d \quad (\textrm{i.e. } d_\Delta \geq 0),
    \label{eqn719}
\end{align}
where $\alpha_{\Delta}, \beta_{\Delta}, d_{\Delta}$ are as in \cref{C7_6_2} (in particular, $\beta_i \not = 0$ for all $i \in \Delta$) and $V_{\Delta}$ is as in \cref{C7_4_1}.
 
\begin{definition}[Minimal $\Delta$]
    \label[definition]{definition7718} 
    Let $\Delta$ satisfy conditions (\ref{eqn718}) and (\ref{eqn719}) above. Define $V_{\Delta_{\alpha}} : = \{ i: i \in  V_{\Delta} \textrm{ and } \alpha_i > 0\} = \{ i:\alpha_i > 0\}, V_{\Delta_{\beta}}:=\{i: i \in V_{\Delta} \textrm{ and } \beta_i > 0\}$ so $V_{\Delta} = V_{\Delta_{\alpha}} \cup V_{\Delta_{\beta}}$. We say $\Delta$ is \textit{minimal} in $S_f(0)$ if for all $i \in V_{\Delta_{\beta}}$ 
    \begin{align}
    (\sum_{j \in V_{\Delta_{\beta}} \backslash \{i\}} \beta_j) - \alpha_{\Delta} < 0.
    \label{eqn720}
    \end{align}
\end{definition}

It follows that
$$
S_f(0) = \bigcup Z_{\Delta},
$$
where the union is over the minimal $\Delta$. Note $Z_{\Delta}$ are the components of $S_f(0)$.

\begin{definition}
    Let $\Delta$ be an arbitrary cone (face) of $\Sigma$ and define
    \begin{align*}
        \Gamma_{\Delta} : = min\{\alpha_{\Delta}, \beta_{\Delta}    \},
        \\\Omega_{\Delta} :=max\{  \alpha_{\Delta}, \beta_{\Delta}  \}.
    \end{align*}
    In particular, for the cone $\Sigma$ defining $M$, $(\Gamma_{\Sigma}, \Omega_{\Sigma}) =(|\alpha|, |\beta|)$.
\end{definition}

\begin{lemma}
    \label[lemma]{lemma78113} 
    $S_f(0) \subset Sing(X)$, where $Sing(X)$ is the locus of non-regular points of $X$. Moreover, $S_f(0) \subset Sing_{\mathbb{Z}}(X)$, where $Sing_{\mathbb{Z}}(X)$ is the locus of points that are not smooth in their fibers.
\end{lemma}

\textit{Proof.} Let $x \in S_f(0) \subset X \subset M$ and $\Delta$ be minimal with $x \in Z_{\Delta}$. Let $\mathfrak{p} \subset \mathbb{Z}[x_1, \dots, x_r, {y_{r+1}}^{\pm}, \dots, {y_n}^{\pm}]$ be the corresponding prime ideal. Then by (\ref{eqn719}), $\mathfrak{p}$ contains $u_1, \dots, u_k, v_{i_1}, \dots, v_{i_l}$ such that $\sum_{j=1}^l \beta_{i_j} \geq d$. It follows that as the $u,v$-variables form part of a regular sequence of parameters in $\mathcal{O}_{M, x}$, $ord_x(f) \geq d \geq 2$. Thus $\mathcal{O}_{X,x} \cong \mathcal{O}_{M,x}/(f)$ cannot be a regular local ring as if it were then $(f)$ could be generated by regular parameters in $ \mathcal{O}_{M,x} $ (see \cite[Lemma 00NR]{stacks-project}) and $ord_x(f)$=1. Now suppose $x$ corresponds to $x_p \in X_p$. The $u_1, \dots, u_k, v_{i_1}, \dots, v_{i_l}$ again form part of a regular sequence in $\mathcal{O}_{M_p, x_p}$ so $ord_{x_p}(f_p) \geq d \geq 2$. Thus by the same arguments $\mathcal{O}_{X_p,x_p} \cong \mathcal{O}_{M_p,x_p}/(f_p)$ is not regular. As $X_p$ is a variety over a perfect field this is equivalent to $x_p$ not being smooth and so $x \in Sing_{\mathbb{Z}}(X)$. \qedsymbol


\subsubsection{Blowing up and strict transform}
Let $Z_{\Delta} \subset S_{f}(0)$ so that $\Delta$ satisfies (\ref{eqn718}), (\ref{eqn719}) and let $\pi: B_{Z_{\Delta}} M \rightarrow M$ be the blow-up. As in \cref{C7_6_2} let $X'$ be the strict transform and $U_{\sigma_i}$ an affine chart of $B_{Z_{\Delta}} M$ corresponding to an $i \in V_{\Delta}$ and maximal cone $\sigma_i$. Then $X' \cap U_{\sigma_i}$ is given by $(w^{\alpha''}-w^{\beta''}y^{\gamma})$, where $\alpha_j'' = \alpha_j, \beta''_j= \beta_j$ for $ j \not = i$ and  $\alpha_i'' = 0, \beta_i'' = d_{\Delta}$. It follows that
$$
|\alpha''| = 
\begin{cases}
  d & \textrm{if }i \in V_{\Delta_{\beta}}\\
  d-\alpha_i & \textrm{if } i \in V_{\Delta_{\alpha}}
\end{cases},
\quad 
|\beta''| = 
\begin{cases}
  |\beta|+d_{\Delta}-\beta_i & \textrm{if } i \in V_{\Delta_{\beta}}\\
  \beta_{\Delta}+d_{\Delta} & \textrm{if } i \in V_{\Delta_{\alpha}}
\end{cases}.
$$

If $\Delta$ is minimal then with respect to the lexicographic ordering $(|\alpha''|, |\beta''|) < (|\alpha|=d, |\beta|)$ and so for each cone $\sigma \in \Sigma'$, $
(\Gamma_{\sigma}, \Omega_{\sigma}) < (\Gamma_{\Sigma}, \Omega_{\Sigma}).$


\subsection{General case}
We extend the constructions in the previous subsection and prove \cref{theorem7201} when $X$ is a hypersurface.


\subsubsection{Maximal order}
Let $\sigma \in \Sigma$ be a cone and let $\tau \in \Sigma$ be a cone with $\sigma$ as a face. Then $U_{\tau}$ is an affine open subscheme of $M$ and $U_{\tau} \cap X$ is given by a binomial $u^{\alpha}-v^{\beta}y^{\gamma}$ as in the affine case. $\sigma$ will correspond to a $V_{\sigma}$ in this chart and define as before
\begin{align}
    &\alpha_{\sigma}: =\sum_{i \in V_{\sigma}} \alpha_i, \quad \beta_{\sigma}:= \sum_{i \in V_{\sigma}} \beta_i
    \\&\Gamma_{\sigma} : = min\{\alpha_{\sigma}, \beta_{\sigma}    \}, \quad \Omega_{\sigma} :=max\{  \alpha_{\sigma}, \beta_{\sigma}  \}.
    \label{eqn7.22} 
\end{align}

Define
\begin{align}
 &\Gamma_{\Sigma}: = \underset{\sigma \in \Sigma}{\textrm{max }}\Gamma_{\sigma}, \quad  V_{\Sigma}:=\{\sigma \in \Sigma: \Gamma_{\sigma}=\Gamma_{\Sigma}\}.
\end{align}

Note $\Gamma_{\Sigma} = 0$ implies $X$ is smooth over $\textrm{Spec} \, \mathbb{Z}$. 

\begin{definition}
    Let $\Delta \in \Sigma$ be a cone. We say $Z_{\Delta}$ (or $\Delta)$ is \textit{admissible} if $\Delta \in V_{\Sigma}$, i.e.  $\Gamma_{\Delta}=\Gamma_\Sigma$. We say $\Delta$ is \textit{minimal} if it is admissible and $\Gamma_{\Delta_1} < \Gamma_{\Delta}$ for every proper face $ \Delta_1$ of $\Delta$.
\end{definition}

Let $\sigma \in \Sigma$ with $\Delta$ a face. Then in $U_{\sigma}$, $X \cap U_{\sigma}$ is given by a $u^{\alpha}-v^{\beta}y^{\gamma}$ with $\Gamma_{\sigma}=|\alpha| \leq |\beta| = \Omega_{\sigma}$. If $\Delta \in \Sigma$ is minimal then locally at $U_{\sigma}$, $\Delta$ is minimal in the sense of \cref{definition7718}, i.e. conditions (\ref{eqn718}), (\ref{eqn719}), (\ref{eqn720}) are satisfied. Conversely, if $\sigma \in \Sigma$ is admissible and $\Delta$ a face that is minimal in the sense of \cref{definition7718}, then $\Delta$ is globally minimal.
\\
\\Define
\begin{align}
    &V_{min, \Sigma} := \{\sigma \in \Sigma: \quad \sigma \textrm{ minimal } \}
    \label{eqn726}    
    \\&\Omega_{\Sigma}:=\underset{\sigma \in V_{min, \Sigma}}{\textrm{max}}{\Omega}_{\sigma}, \quad W_{\Sigma}:=\{\sigma \in V_{min, \Sigma}: {\Omega}_{\sigma}={\Omega}_{\Sigma}\}.
    \label{eqn727} 
\end{align}



\subsubsection{Ideal sheaf of Hasse derivatives}
\label{C7_8_2_2}  
We glue together the local closed subschemes defined by Hasse derivative monomials into a global closed subscheme of $M$. The components of the closed subscheme will be the blow-up centres in our resolution.
\\
\\For each $\sigma \in V_{\Sigma}$ we have the affine $U_{\sigma}$ where $X \cap U_{\sigma}$ is given by the binomial $f_{\sigma}$ with $|\alpha| = d: = \Gamma_{\Sigma}$. Write $I_{\sigma}$ for the ideal in $U_{\sigma}$ given by the Hasse derivative monomials of $f_{\sigma}$ of order $< d$, which defines $S_{f_{\sigma}}(0)$. On $U_{\sigma}$ with $\sigma$ not admissible we define $I_{\sigma}$ to be the unit ideal. 

\begin{lemma}
    \label[lemma]{lemma7825} 
    The $I_{\sigma}$ glue along faces.
\end{lemma}

\textit{Proof.} Let $\sigma \in \Sigma$ and $\Delta \in \Sigma$ be a face of $\sigma$ so $U_{\Delta}$ is given by localizing in $U_{\sigma}$ each $x_i, i \not \in V_{\Delta}$, and in $U_{\sigma}$, $\Gamma_{\Delta} \leq d$. When $\sigma$ is not admissible $\Delta$ cannot be admissible. Assume $\sigma$ is admissible and suppose $\Delta$ is not admissible so $\Gamma_{\Delta} < d$. Without loss of generality assume $\Gamma_{\Delta}=\alpha_{\Delta}$ and consider 
$$
\zeta = (\zeta_1, \dots, \zeta_k), \quad \zeta_i = 
\begin{cases}
    \alpha_i & \textrm{if } i \in V_{\Delta_{\alpha}}\\
    0 & otherwise
\end{cases}.
$$
$u^{\alpha-\zeta}$ will be a Hasse derivative monomial of order $< d$ and hence an element of $I_{\sigma}$. But clearly $u^{\alpha-\zeta}$  will become a unit in $U_{\Delta}$ and so the localization of $I_{\sigma}$ defines the whole ring. If $\Delta$ is admissible then $\alpha_{\Delta}=d$ and only $v$-variables are inverted. Thus the Hasse derivative monomials of order $< d$ in $U_{\Delta}$ in the $u$-variables are the same as those in $U_{\sigma}$. Now $v^{\beta}y^{\gamma}$ in the binomial $f_{\sigma}$ becomes $v^{\beta'}y^{\gamma'}$ in the binomial $f_{\Delta}$, where $\beta'$ can be thought of as $\beta$ but removing $\beta_i: i \not \in V_{\Delta}$ and $\gamma'$ can be thought of as $\gamma$ but including $\beta_i:i \not \in V_{\Delta}$ as these variables will now be included among the invertible $y$'s. Any Hasse derivative monomial $v^{\beta'-\delta}$ of order $|\delta| < d$ will correspond to the localization of the Hasse derivative monomial $v^{\beta-\delta}$, where we view $\delta $ as an element of $ \mathbb{N}^{r-k}$ by appending $0$'s. Conversely, any Hasse derivative monomial $v^{\beta-\delta}$ of order $|\delta| < d$ under localization becomes a Hasse derivative monomial $v^{\beta'-\delta'}$ of order $|\delta'| \leq |\delta| < d$ and this is never a unit as $|\delta| < d \leq |\beta|$. \qedsymbol

\begin{definition}
    Define $S$ to be the closed subscheme of $M$ defined by the gluing of the $S_{f_{\sigma}}(0)$ for $\sigma$ admissible as in the lemma above. 
\end{definition}

\begin{observation}
    \label{observation7827} 
    By \cref{lemma78113}, $S \subset Sing_{\mathbb{Z}}(X) \subset X$. Moreover, as $S$ is obtained by gluing along faces, which is a gluing along underlying monoid subschemes, $S$ is the $\mathbb{Z}$-realization of an equivariant closed monoid subscheme of $M_{Sch}$ and $X_{Sch}$. In particular, it is a $\mathbb{Z}$-flat $\Lambda$-subscheme.
\end{observation}



\subsubsection{Proof of \cref{theorem7201} for hypersurfaces}
The proof will be as in \cite[Section 5.2]{BM1}. Unlike \cite[Section 5.2]{BM1} we use the collection of cones $W_{\Sigma}$ (see (\ref{eqn727})) and an explicit lexicographic triple that decreases after blowing up. Let $\Delta \in W_{\Sigma}$ and consider the blow-up of $M$ along $Z_{\Delta}$ (a component of $S$), which corresponds to a refinement $\Sigma'$. In particular, $\Delta$ is minimal so it follows that
$$
(\Gamma_{\Sigma'}, \Omega_{\Sigma'}, |W_{\Sigma'}|) < (\Gamma_\Sigma, \Omega_{\Sigma}, |W_{\Sigma}|).
$$
Admissible (and hence minimal) $Z_{\Delta}$ are contained in $S \subset X$ and so each $Z_{\Delta} \cap X = Z_{\Delta}$ is already smooth over $\textrm{Spec} \, \mathbb{Z}$. Thus by observation \ref{observation7827} and \cref{C7_7} all conditions of \cref{theorem7201} on the blow-up centres are satisfied by minimal $Z_{\Delta}$.
It remains to show we can make $X$ smooth. As $|W_{\Sigma}|$ is finite, after finitely many blow-ups by minimal $Z_{\Delta}$ there exists a refinement $\Sigma_m$ with $(\Gamma_{\Sigma_m}, \Omega_{\Sigma_m}) < (\Gamma_{\Sigma_m}, \Omega_{\Sigma_m})$. After finitely many more such blow-ups we have a $\Sigma_n$, $X_n$, $M_n$ such that $\Gamma_{\Sigma_n} = 0$ and so $X_n$ is smooth over $\textrm{Spec} \,  \mathbb{Z}$. \qedsymbol


\section{Order of binomials and standard bases}
\label{C7_9} 
In this section we prove facts about the binomial equations defining the embedding of the main theorem and pullbacks to fibers. We relate the order of a binomial at a point with the point in its fiber as in \cite[Section 3.1]{BM1} and introduce a basis of binomials, which have as pullbacks the standard bases of \cite[Section 6]{BM1}. 

\begin{lemma} 
    \label[lemma]{lemma7912} 
    Let $X = \textrm{Spec} \, \mathbb{Z}[x_1, \dots, x_r, y_{1}^{\pm}, \dots, y_{n-r}^{\pm}]$ and $a \in X$ corresponding to $a_p \in X_p$. If $f \in \mathbb{Z}[x_1, \dots, x_r, y_{1}^{\pm}, \dots, y_{n-r}^{\pm}]$ then
    $$
    ord_a(f) \leq ord_{a_p}(f_p) \leq ord_{x(a_p)}(f_p),
    $$
    where $ord_a(f)$ is the usual order, $ord_{a_p}(f_p)$ is the order at $a_p$ of the pullback in the fiber over $(p)$ and $ord_{x(a_p)}(f_p)$ is the alternate notion of order described in \cite[Section 3.1]{BM1} computed on the fiber $X_p$. Moreover, if $f$ is a monomial we get equality.
\end{lemma}

\textit{Proof.} Follows from basic facts on order and fibers and \cite[Section 3.1]{BM1}. \qedsymbol
\\
\\Let $f=x^{\alpha}-x^{\beta}y^{\gamma} \in \mathbb{Z}[x_1, \dots, x_r, y_{1}^{\pm}, \dots, y_{n-r}^{\pm}], |\alpha| \leq |\beta|$ be a binomial where $x^{\alpha}, x^{\beta}$ share no common variables so we can write $f=u^{\alpha}-v^{\beta}y^{\gamma}$ as before.  It follows from \cite[Corollary 3.5]{BM1} and \cref{lemma7912} that
\begin{align*}
    S_f(0):&=\{a \in \textrm{Spec} \, \mathbb{Z}[x_1, \dots, x_r, y_{1}^{\pm}, \dots, y_{n-r}^{\pm}]: ord_a(f) = |\alpha|\}\\
    &= \{ a \in \textrm{Spec} \, \mathbb{Z}[x_1, \dots, x_r, y_{1}^{\pm}, \dots, y_{n-r}^{\pm}]: \space \textrm{all } u_i \textrm{ vanish on } a \textrm{ and} \sum_{i:v_i \in a}\beta_i \geq |\alpha|\}.
\end{align*}

Thus, we may decompose
$$
S_f(0) = \bigsqcup_{(p) \in \textrm{Spec} \, \mathbb{Z}} (S_{f_p}(0) \subset X_p),
$$
where $X=\textrm{Spec} \, \mathbb{Z}[x_1, \dots, x_r, y_{1}^{\pm}, \dots, y_{n-r}^{\pm}]$ and $S_{f_p}(0)$ is as in \cite[Section 3.2]{BM1}. 

\begin{remark}
    It follows that $S_f(0)$ is the same underlying set as that defined by Hasse derivatives in \cref{C7_8_1_2} (cf. \cite[Remark 3.7]{BM1}). Further, the global closed subscheme $S \subset M$ defined by Hasse derivatives in \cref{C7_8_2_2} can be seen as the locus of points of maximal order $\Gamma_{\Sigma}$. This locus will similarly decompose into the locus in each fiber (cf. \cite[Section 5.2]{BM1}). This means that one could rewrite \cref{C7_8} using components of equimultiple loci as in \cite{BM1}.
\end{remark}    


Let $X \hookrightarrow M$ be the closed embedding of toric $\Lambda$-schemes of \cref{theorem7201} in the affine case as in \cref{C7_5}. Let $I$ be the ideal defining $X$ and recall that $I$ is generated by binomials of the form $x^{\alpha}-x^{\beta}y^{\gamma}, 1-y^{\gamma}, \alpha,\beta \in \mathbb{N}^{n-m}, \gamma \in \mathbb{Z}^{m}, 1 \leq |\alpha| \leq |\beta|$. Let $J \subset \mathbb{Z}[{y_1}^{\pm},\dots, {y_{n-r}}^{\pm}]$ denote the ideal generated by all $1-y^{\gamma^j} \in I$. We show that $I$ is generated mod $J$ by certain binomials, the \textit{standard basis elements}. We will use the notation of \cite[Section 6]{BM1} and remark on how the arguments in this section apply to toric $\Lambda$-schemes. Further, writing $I_p$ for the pullback of $I$ defining the fiber $X_p$, the standard basis elements of $I$ will be the standard basis (\cite[Remark 6.6]{BM1}) of $I_p$ for all prime ideals $(p)$. Recall
$$
X \cong \textrm{Spec}\left( \frac{\mathbb{Z}[x_1, \dots, x_{r}, {y_1}^{\pm},\dots, {y_{n-r}}^{\pm}]}{I=(x^{\alpha^i}-x^{\beta^i}y^{\gamma^i}, 1-y^{\gamma^j} , i=1, \dots, r, j=1, \dots, s)} \right)
$$
for some $r,s$.
\\
\\\textbf{Distinguished point.} Let $X$ be as above. Let $a$ be the point in $\mathbb{A}^n_{\mathbb{Z}}$ corresponding to $\mathfrak{p}=(x_1, \dots, x_{r}, y_1-1, \dots, y_{n-r}-1)$. Then $a \in X$ and as $\mathfrak{p}$ contains no primes $p$, $a$ corresponds to $a_0 \in X_0$, the generic fiber. Now write $(p, a) = (p, x_1, \dots, x_{r}, y_1-1, \dots, y_{n-r}-1)$ for any prime ideal $(p) \in \textrm{Spec} \, \mathbb{Z}$ ($p=0$ gives $a$). Then $(p, a) \in X$ (it is a closed point if $p \not = 0$) and corresponds to $(p,a)_p \in X_p$, which will be the distinguished point of $X_{p}$ (see \cite[Section 6.2]{BM1}). Note that $a$ will be a $\mathbb{F}_1$-point (see \cite{borger2009lambdarings}) and does not correspond to a prime ideal of the underlying monoid.

\begin{example}[Intersection distinguished]
     \label[example]{example71001} 
     Let $X$ be a smooth normal toric $\Lambda$-scheme with fan $\Sigma$ and $E$ be the locally toric simple normal crossings divisor of \cref{example7412}. On any affine open $U_{\sigma}$ the distinguished point $a$ will be contained in every codimension one orbit closure intersecting $U_{\sigma}$. Then $E$ is intersection distinguished for the $\Lambda$-cover given by cones (see \cite[Definition 2.3.19]{OrderReductionLambdamarkedMonomialIdeals}).
\end{example}

\textbf{Diagram of initial vertices.} Let $A$ be a commutative ring and $A\llbracket  X \rrbracket=A\llbracket  X_1, \dots, X_n \rrbracket$ be the ring of formal power series. We use the total order on $\mathbb{N}^q$ that is given by the lexicographic ordering of $(q+1)$-tuples $(|\alpha|, \alpha_1, \dots, \alpha_q)$ for $\alpha=(\alpha_1, \dots, \alpha_q) \in \mathbb{N}^q$. Let $F= \sum_{\alpha \in \mathbb{N}^q}F_{\alpha}X^{\alpha} \in A\llbracket X \rrbracket$ and define supp$F := \{\alpha: F_{\alpha} \not = 0   \}$. The \textit{initial exponent}, exp $F$, is the smallest element of supp$F$ given by the lexicographic ordering ($exp F : = \infty$ if $F=0$). Let $I \subset A\llbracket X \rrbracket$ be an ideal. The \textit{diagram of initial exponents} $\mathcal{R}(I) \in \mathbb{N}^q$ is the subset
$$
\mathcal{R}(I) : = \{ \textrm{exp}F: F \in I \backslash \{0\} \}.
$$
It follows that $\mathcal{R}(I)+\mathbb{N}^q=\mathcal{R}(I)$ and so there exists a smallest finite subset $\mathcal{B}$ of $\mathcal{R}(I)$, the \textit{vertices} of $I$, such that $\mathcal{R}(I)=\mathcal{B}+\mathbb{N}^q$.

\begin{lemma}[Standard bases]
    \label[lemma]{lemma7902} 
    $I$ is generated modulo $J \cdot \mathbb{Z}[{y_1}^{\pm},\dots, {y_{n-r}}^{\pm}][x_1, \dots, x_{r}]$ by a set of binomials $F^i = x^{\alpha^i}-x^{\beta^i}y^{\gamma^i}$. The $\{F^i\}$ will be called the \textit{standard basis} of $I$ mod $J \cdot \mathbb{Z}[{y_1}^{\pm},\dots, {y_{n-r}}^{\pm}][x_1, \dots, x_{r}]$. For each prime ideal $(p) \subset \mathbb{Z}$ the standard basis of $I$ mod $J \cdot \mathbb{Z}[{y_1}^{\pm},\dots, {y_{n-r}}^{\pm}][x_1, \dots, x_{r}]$ is (mod $p$) the standard basis (\cite[Remark 6.6]{BM1}) of $I_p$ mod $J_p \cdot \mathbb{F}_p[{y_1}^{\pm},\dots, {y_{n-r}}^{\pm}][x_1, \dots, x_{r}]$.
\end{lemma}

\textit{Proof.} The constructions in \cite[Section 6]{BM1} over perfect fields apply naturally to the binomials in $\mathbb{Z}[{y_1}^{\pm},\dots, {y_{n-r}}^{\pm}][x_1, \dots, x_{r}]$. Note these constructions rely only on the data of exponents of binomials. Hironaka division holds over $\mathbb{Z}$ by considering the embedding $\mathbb{Z}[{y_1}^{\pm},\dots, {y_{n-r}}^{\pm}][x_1, \dots, x_{r}] \hookrightarrow \mathbb{Q}[{y_1}^{\pm},\dots, {y_{n-r}}^{\pm}][x_1, \dots, x_{r}]$ and it follows we can produce a set of binomials with analogous properties to the standard basis of \cite[Remark 6.6]{BM1}. As the exponents of binomials and the diagram of initial vertices are unchanged under the pullback to fibers, (mod $p$) we get the standard basis of \cite[Remark 6.6]{BM1}. \qedsymbol


\section{Hilbert-Samuel functions and strata}
\label{C7_11} 
We prove properties of Hilbert-Samuel functions and maximal strata analogous to \cite[Section 7]{BM1}. We describe how these relate to fibers, which we will make use of in the proof of \cref{theorem7201} in the following section.


\subsection{Hilbert-Samuel functions and fibers}
Let $R$ denote a Noetherian local ring with maximal ideal $\mathfrak{m}$. The \textit{Hilbert-Samuel function} $H_R: \mathbb{N} \rightarrow \mathbb{N}$ is the function defined by
$$
H_R(l) = \textrm{length} \frac{R}{\mathfrak{m}^{l+1}}, \quad l \in \mathbb{N}.
$$
We partially order the set of functions $N:=\{ H:\mathbb{N} \rightarrow \mathbb{N} \}$ by their values, i.e. if $H, H' \in N$ then $H \leq H'$ if and only if $H(l) \leq H'(l)$, for all $l \in \mathbb{N}$. We define the Hilbert-Samuel function $H_{X, b}$ of a Noetherian locally ringed space $X=(|X|, \mathcal{O}_X)$ at a point $b$ as the Hilbert-Samuel function of the local ring $\mathcal{O}_{X, b}$. If $k \cong R/\mathfrak{m}$, then $H_R(l)=dim_{k}R/\mathfrak{m}^{l+1}$ for all $l$. In particular, if $I \subset k\llbracket x \rrbracket=k\llbracket x_1, \dots, x_n \rrbracket$ is an ideal in the ring of formal power series, then $k\llbracket x \rrbracket/I$ is a Noetherian local ring with maximal ideal $\mathfrak{m}+I$, where $\mathfrak{m}$ is generated by the variables. Thus, it follows by Hironaka division in the context of formal power series \cite[Theorem 3.17]{https://doi.org/10.48550/arxiv.alg-geom/9508005} that
\begin{align}
   H_{k\llbracket x \rrbracket/I}(l) = \# \{\alpha \in \mathbb{N}^n : \alpha \not \in \mathcal{R}(I), |\alpha| \leq l \}. 
   \label{eqn728}
\end{align}
\\Recall the notation of \cref{C7_9}. Let $X \hookrightarrow M$ be a closed embedding of affine toric $\Lambda$-schemes as before and let $I \subset \mathbb{Z}[x_1, \dots, x_{r}, y_1^{\pm}, \dots, y_{n-r}^{\pm}]$ denote the ideal underlying $X$. Let $J \subset \mathbb{Z}[y_1^{\pm}, \dots, y_{n-r}^{\pm}]$ be the ideal generated by the $1-y^{\gamma}$ and $a$ be the distinguished point of $X$. Let
$$
F^i = x^{\alpha^i}-x^{\beta^i}y^{\gamma^i} \in I, \quad i=1, \dots, s,
$$
be the standard basis elements of $I$ mod $J \cdot \mathbb{Z}[x_1, \dots, x_{r}, y_1^{\pm}, \dots, y_{n-r}^{\pm}]$.

\begin{proposition} 
    \label[proposition]{proposition71111}
    $H_{X,a} = H_{X_p, (p,a)_p}$ for all prime ideals $(p)$, i.e. the Hlbert-Samuel function at the distinguished point of $X$ is the same as the Hilbert-Samuel function at the distinguished points of the toric varieties $X_p$.
\end{proposition}

\textit{Proof.} As $\mathcal{O}_{X, a} \cong \mathcal{O}_{X_0, a_0}$ we have $H_{X,a} = H_{X_0, a_0}$. We will show that for any primes $p,p'$, $H_{X_p, (p,a)_p}=H_{X_{p'}, (p',a)_{p'}}$, which will prove the proposition by letting $p=0$. 
\\Let $X_p$ be the affine toric variety over the perfect field $\mathbb{F}_p=k$ embedded as a closed subscheme in $\mathbb{A}^n_k$. We can compute the Hilbert-Samuel function at the distinguished point by first completing at the stalk then using the isomorphism of \cite[Theorem 3.1]{BM1} to show that it is equal to the Hilbert-Samuel function of the ring 
$$
\frac{k\llbracket x_1 \dots, x_{r}, y_1, \dots, y_{n-r} \rrbracket}{\hat{I}_p}
$$ 
and so may apply (\ref{eqn728}). Here $I_p$ is the pullback of the ideal $I$ defining the embedding of $X$ in $\mathbb{A}_{\mathbb{Z}}^n$ (see observation \ref{observation7522}) and $\hat{I}_p$ is the ideal in the completion generated by $I_p$. Let $J_p$ be the pullback of $J$ and denote by $V(J_p) \subset \textrm{Spec} \, k[{y_1}^{\pm},\dots, {y_{n-r}}^{\pm}]$ the closed subscheme defined by $J_p$. By \cite[Remark 6.7]{BM1} the vertices of $\mathcal{R}(\hat{I}_p)$ will be given by $(\alpha^i, 0) \in \mathbb{N}^{r} \times\mathbb{N}^{n-r}, i = 1, \dots, s$ (the exponents of the standard basis elements of $I_p$ mod $J_p$) and $\nu_p$ elements $(0, {\gamma}^i) \in \mathbb{N}^{r} \times \mathbb{N}^{n-r}$ with $|\gamma^i|=1$, where $\nu_p$ is the codimension of $V(J_p) $ in $\mathbb{A}_k^{n-r}$ at the distinguished point $(a,p)_p$. Note that while $\alpha^i$-vertices are independent of $p$ by \cref{lemma7902} the $(0,\gamma^i)$ come from certain linear terms in the images of $1-y^\gamma$ in the formal power series, which will differ depending on char $k$. Consider the set
$$
R_{l ,p}:=\{\alpha \in \mathbb{N}^n : \alpha \not \in \mathcal{R}(\hat{I}_p), |\alpha| \leq l \} \subset \mathbb{N}^n,
$$
where $\hat{I}_p$ is the ideal of the formal power series ring defined above and the number of elements in this set gives $H_{X_{p}, (p, a)_p}(l)$ by (\ref{eqn728}). Note $\alpha =(\alpha_x, \alpha_y) \in \mathbb{N}^{r} \times \mathbb{N}^{n-r}$ is contained in $R_{l,p}$ if and only if $\alpha_x$ is not contained in the vertices of $\mathbb{N}^{r}$ given by the $\alpha^i$ and $\alpha_y$ has a $0$ in each coordinate given by a $\gamma^j$. Let $p,p'$ be distinct primes. Then by the lemma below $\nu_p = \nu_p'$ and we can define the set map
$$
f: R_{l ,p} \rightarrow R_{l ,p'}, \quad \alpha =(\alpha_x, \alpha_y)  \mapsto \alpha'=(\alpha_x', \alpha_y'),
$$
where $\alpha_x'=\alpha_x$ and $\alpha_y'$ is given by mapping from left to right the non-zero coordinates of $\alpha_y$ (of which there are $n-r-\nu_p$) to (from left to right) the non-zero free coordinates specified by the $\gamma^j$ over $p'$ (of which there are also $n-r-\nu_{p'}$), which implies $\alpha' \in R_{l, p'}$. Clearly $f$ is a bijection and so $H_{X_{p}, (p, a)_p}(l) = H_{X_{p'}, (p', a)_p'}(l)$ for all $l$, which implies $H_{X_{p}, (p, a)_p} = H_{X_{p'}, (p', a)_p'}$. \qedsymbol

\begin{lemma} 
    \label[lemma]{lemma71112}
    With notation as in the preceding proposition $\nu_p = \nu_{p'}$.
\end{lemma}

\textit{Proof.} Fix a $p$. Then 
$$
V(J_p) \cong \textrm{Spec} \, \frac{k[y_1^{\pm}, \dots, y_{n-r}^{\pm}]}{(1-y^\gamma \in I_p)}
$$
and $I_p$ is a prime ideal so $J_p$ is prime in $k[y_1^{\pm}, \dots, y_{n-r}^{\pm}]$. It follows that $V(J_p)$ itself is a toric variety and so corresponds to a monoid algebra of a monoid $A$. In particular, $A$ will be the image of the surjective monoid morphism from $\mathbb{Z}^{n-r}$ induced by the congruence generated by the $1-y^\gamma$. Then $A$ is itself isomorphic to some $\mathbb{Z}^c$ and hence $V(J_p)$ is itself a torus. Further, as this was induced by the underlying monoids this $c$ is independent of $p$. But the codimension of a torus at a point will be the same regardless of the point and equal to the rank of the complement of $\mathbb{Z}^c$ in $\mathbb{Z}^{n-r}$ and hence $\nu_p=\nu_{p'}$ for any other $p'$. \qedsymbol 
\\
\\The Hilbert-Samuel function at the distinguished point $a$ does not quite behave the same way as for toric varieties over perfect fields.

\begin{example} 
    \label[example]{example71113}
    Let $\mathfrak{p}=(p, x_1, \dots, x_n) \subset \mathbb{Z}[x_1, \dots, x_n]$ be the prime ideal corresponding to $(a,p)$ (when there are no $y$-variables) for some $p \not = 0$ and let $\mathfrak{m} = (x_1, \dots, x_n) \subset \mathbb{F}_p[x_1, \dots, x_n]$ be the prime ideal correpsonding to $(p, a)_p$. Fix a $k \in \mathbb{N}$ and consider the Noetherian local rings $R_{\mathfrak{p}}=\mathbb{Z}[x_1, \dots, x_n]_{\mathfrak{p}}, R_{\mathfrak{m}}=\mathbb{F}_p[x_1, \dots, x_n]_{\mathfrak{m}}$ and the surjective morphism induced by the surjective local ring morphism
    $$
    i_p: \frac{R_{\mathfrak{p}}}{(\mathfrak{p} \cdot R_{\mathfrak{p}})^{k+1}} \rightarrow \frac{R_{\mathfrak{m}}}{(\mathfrak{m} \cdot R_{\mathfrak{m}})^{k+1}}.
    $$
    Then one can show any chain of length $l$ in $\frac{R_{\mathfrak{m}}}{(\mathfrak{m} \cdot R_{\mathfrak{m}})^{k+1}}$
    $$
    0 \subsetneq (h_1, \dots, h_r) \subsetneq \dots \subsetneq (h_1, \dots, h_s)
    $$
    lifts to a chain of length $l$ in $\frac{R_{\mathfrak{p}}}{(\mathfrak{p} \cdot R_{\mathfrak{p}})^{k+1}}$. In particular, if we consider the maximal chain we can lift and construct the augmented chain in $\frac{R_{\mathfrak{p}}}{(\mathfrak{p} \cdot R_{\mathfrak{p}})^{k+1}}$
    $$
    0 \subsetneq (p^{k}) \subsetneq (p^{k-1}) \subsetneq \dots \subsetneq (p) \subsetneq (p, h_1) \subsetneq \dots, \subsetneq (p, h_1, \dots, h_s)
    $$
    of length $k$ greater. Thus we have shown $H_{\mathbb{A}^n_{\mathbb{Z}}, (p, a)}(k) \geq H_{\mathbb{A}_{\mathbb{F}_p}^n, (p, a)_p}(k)+k$ and so $H_{\mathbb{A}^n_{\mathbb{Z}}, (p, a)} > H_{\mathbb{A}_{\mathbb{F}_p}^n, (p, a)_p}$.
\end{example}

\begin{corollary} 
    \label[corollary]{corollary71114} 
    Combined with \cref{proposition71111} we get $H_{\mathbb{A}^n_{\mathbb{Z}}, (p, a)} > H_{\mathbb{A}^n_{\mathbb{Z}}, a} $. Thus it is not true that for all $b \in \mathbb{A}^n_{\mathbb{Z}}$, $H_{\mathbb{A}^n_{\mathbb{Z}}, b} \leq H_{\mathbb{A}^n_{\mathbb{Z}}, a}$ as opposed to  in any fiber by \cite[Theorem 7.1 (1)]{BM1}.
\end{corollary}


\subsection{Maximal strata}
\label{C7_11_2} 
Let $X \hookrightarrow M$ be a closed embedding of affine toric $\Lambda$-schemes with 
$$
I \subset \mathbb{Z}[x_1, \dots, x_{r}, y_1^{\pm}, \dots, y_{n-r}^{\pm}] \textrm{ and } J \subset \mathbb{Z}[y_1^{\pm}, \dots, y_{n-r}^{\pm}]
$$ 
as before. Let
$$
F^i = x^{\alpha^i}-x^{\beta^i}y^{\gamma^i} \in I, \quad i=1, \dots, s
$$
be the standard basis elements of $I$ mod $J \cdot \mathbb{Z}[x_1, \dots, x_{r}, y_1^{\pm}, \dots, y_{n-r}^{\pm}]$. For each $i$, if $|\alpha^i| = 1$ then $x^{\alpha^i}=x_{j(i)}$ for some $j(i) \in \{1, \dots, r\}$, and $x_{j(i)}$ occurs (to nonzero power) in no monomial $x^{\alpha^j}, j \not = i$, and in no monomial $x^{\beta^j}$.
\\Define $N=N(I)$ to be the closed-subscheme of $M$ defined by the binomials $1-y^\gamma \in J$ together with $F^i$ for all $i=1, \dots, s$ such that $|\alpha^i|=1$. Then $N$ will be the $\mathbb{Z}$-realization of the closed monoid subscheme of $M_{MSch}$ defined by the congruence generated by the binomials and so will be smooth over $\textrm{Spec} \, \mathbb{Z}$. After reordering the indices $i$ if necessary, we can assume that
\begin{enumerate}
    \item $|\alpha^i| \geq 2, i = 1, \dots, t$ and $|\alpha^i| = 1$, $i=t+1, \dots, s$, $t \leq s$.
    \item $x_1, \dots, x_k$ are those variables occuring (to nonzero power) in some $x^{\alpha^i}, i=1, \dots, t$, where $k \leq r-(s-t)$. 
\end{enumerate}
We can then view $N$ as a smooth toric $\Lambda$-scheme
$$
N \cong \textrm{Spec} \, \mathbb{Z}[z_1, \dots, z_{r-(s-t)}, w_1^{\pm}, \dots, w_{c}^{\pm}],
$$
where $N \rightarrow M$ is defined by the ring morphism sending $x_i$ to $z_i$ for $ i \not \in \{ t+1, \dots, s \}$ and the description of the torus given by $J$ in the proof of \cref{lemma71112}. In particular, $N$ is a toric $\Lambda$-scheme with $N \hookrightarrow M$ induced by a closed embedding of monoid schemes that is unpointed. Further, by \cref{lemma7902} the fiber $N_p$ for any $(p)$ will be the smooth subvariety constructed from $X_{p}$ in \cite[Section 7.2]{BM1}.
\\
\\Define the \textit{fibered Samuel stratum $S_X(b)$} of $b \in X$ corresponding to $b_p \in X_{p}$ as 
$$
S_X(b) : = \{c \in X \textrm{ that corresponds to a } c_{p'}: H_{X_{p'}, c_p'} = H_{X_{p}, b_p}    \}.
$$

\begin{remark}
    We use fibered Samuel strata instead of the Samuel strata in \cite[Section 7.2]{BM1} as we wish to capture all points that contribute to the maximal Hilbert-Samuel functions in their fibers. By \cref{corollary71114}, if we defined the strata by Hilbert-Samuel functions in $X$, the maximal Hilbert-Samuel function would not necessarily occur at the distinguished point $a$. Then, when we stratified we would not include points such as $(p, a)$, which become the distinguished point in their fibers. 
\end{remark}

\begin{theorem}
    \label[theorem]{theorem71122} 
    (See \cite[Theorem 7.1]{BM1}).
    Let $S_X(a)$ be the fibered Samuel stratum at the distinguished point $a$. Then as a subset of $M$
    $$
    S_X(a) = \{b \in N: x_j \in b, \quad j=1, \dots, k \textrm{ and } \sum_{x_l \in b} (\beta^i)_l \geq |\alpha^i|, \quad i=1, \dots, t\}.
    $$
\end{theorem}

\textit{Proof.} Follows from \cref{proposition71111}, \cref{lemma7912} and \cite[Theorem 7.1 (2)]{BM1}. \qedsymbol

\begin{corollary}
    $S_X(a) \subset X $. Further, $S_X(a) = \bigsqcup_{(p) \in \textrm{Spec} \, \mathbb{Z}} S_{X_{p}}((p,a)_p)$ as sets, where $ S_{X_{p}}((p,a)_p)$ is the usual Samuel stratum in \cite[Section 7.2]{BM1}.
\end{corollary}

\textit{Proof.} The first part follows from \cref{theorem71122} and how $X$ is embedded into $N$. The second part again follows from \cite[Theorem 7.1 (2)]{BM1}. \qedsymbol

\begin{corollary}
    \label[corollary]{corollary71124} 
    $S_X(a) \subset Sing(X) \subset X$,
    where $Sing(X)$ is the locus of non-regular points in $X$. Moreover, $S_X(a) \subset Sing_{\mathbb{Z}}(X)$, where $Sing_{\mathbb{Z}}(X)$ is the locus of points that are not smooth in their fibers.
\end{corollary}

\textit{Proof.} Let $b \in S_X(a) \subset X$. As $N$ is a smooth scheme over $\mathbb{Z}$, $\mathcal{O}_{N, b}$ is a regular local ring. By definition the embedding  $X \hookrightarrow N$ is defined by the ideal $L=(F^i, i=1, \dots, t)$, which implies $\mathcal{O}_{X, b} \cong \mathcal{O}_{N,b}/L_b$. But by \cref{theorem71122} for each $i=1, \dots, t$ the order in $N$, $ord_b(F_i) \geq |\alpha^i| \geq 2$, and so $ord_bL \geq 2.$ As in the proof of \cref{lemma78113} if $\mathcal{O}_{X, b}$ were a regular local ring then $L_b$ could be generated by regular parameters and $ord_b L = 1$. Now suppose $b$ corresponds to $b_p \in X_p$. Then as $X_p \hookrightarrow N_p$ is defined by $L+(p)$ and $ord_{b_p}(L+(p))) \geq ord_b(L) \geq |\alpha^i| \geq 2$, $\mathcal{O}_{X_p, b_p}$ is not a regular ring. It follows that $b_p$ is not a smooth point in $X_p$ and so $b \in Sing_{\mathbb{Z}}(X)$.  \qedsymbol

\begin{corollary}
    \label[corollary]{corollary71126} 
    Let $L \subset \mathbb{Z}[x_1, \dots, x_{r}, y_1^{\pm}, \dots, y_{n-r}^{\pm}]$ be the ideal generated by each $1-y^{\gamma} \in J, F^i, i=t+1, \dots, s$ and for each $i=1, \dots, t$ the Hasse derivative monomials of $F^i$ of order $< |\alpha^i|$. Then as subsets of $M$
    $$
    V(L) = S_X(a)
    $$
    and we can equip $S_X(a)$ with a closed subscheme structure of $M$ or $X$. Moreover, with this scheme structure $S_X(a)$ is the $\mathbb{Z}$-realization of an equivariant closed monoid subscheme of $X_{MSch}$.
\end{corollary}

\textit{Proof.} It suffices to show agreement in $N$, which allows us to argue as in the hypersurface case. Now $N$ is a smooth toric $\Lambda$-scheme and the ideal induced by $L$ defining the closed subscheme in $N$ is generated by the Hasse derivative monomials of the class of each $F^i, i=1, \dots, s$. Thus the agreement of sets follows from \cref{theorem71122} and applying the results of \cref{C7_9} to each $F^i$. As Hasse derivative monomials come from underlying monoid elements the ideal sheaf induced by $L$ in $X$ (or in $M$) comes from the $\mathbb{Z}$-realization of a monoid ideal. The closed subscheme structure is then the $\mathbb{Z}$-realization of an equivariant closed monoid subscheme. \qedsymbol

\begin{lemma} 
    \label[lemma]{lemma71141} 
    Let $N \rightarrow M$ be a closed embedding of smooth toric $\Lambda$-schemes, i.e. with the same conditions as in \cref{theorem7201} but $N$ is also a smooth normal toric $\Lambda$-scheme. Then
    \begin{enumerate}
        \item Every orbit closure of $N$ is the intersection with $N$ of an orbit closure of $M$ that is "transverse" to $N$ (see \cite[Section 2.3]{OrderReductionLambdamarkedMonomialIdeals}).
        \item The intersection with $N$ of an orbit closure of $M$ has irreducible components given by orbit closures of $N$.
    \end{enumerate}
\end{lemma}

\textit{Proof.} Follows from the results of \cref{C7_5} and the local description of $N(I)$ in \cref{C7_11_2}. \qedsymbol


\section{Proof of the main theorem}
\label{C7_12}
Let $X \hookrightarrow M$ be the $\Lambda$-equivariant embedding of toric $\Lambda$-schemes as in \cref{theorem7201}. In this section we prove the theorem following \cite[Section 8]{BM1}. We construct analogous objects in the category of $\Lambda$-schemes and use the results of \cref{C7_11}, \cref{C7_9} and the order reduction of \cite[Section 4]{OrderReductionLambdamarkedMonomialIdeals} to employ the inductive arguments of \cite[Section 8]{BM1} simultaneously in each fiber.


\subsection{\texorpdfstring{$\mathbb{Z}$}{ℤ}-maximal Hilbert-Samuel functions}
Let $H$ be the Hilbert-Samuel function on $X$ maximal amongst those points corresponding to the distinguished point on some affine chart. We call such a $H$ the \textit{$\mathbb{Z}$-maximal} Hilbert-Samuel function of $X$ (cf. \cref{corollary71114}, \cite[Section 8]{BM1}).  

\begin{lemma}
    For all prime ideals $(p) \subset \mathbb{Z}$, $H=H_{p}$, where $H_{p}$ is the maximal Hilbert-Samuel function (taken across all points) for $X_{p}$ (see \cite[Section 8]{BM1}).
\end{lemma}

\textit{Proof.} Let $a_{\sigma}$ be the distinguished point in some affine chart $X_{\sigma} : = X \cap U_{\sigma}$, where $H_{X_{\sigma}, a_{\sigma}} = H_{(X_{\sigma})_{p}, (a_{\sigma}, p)_{p}}$ for all $p$ by \cref{proposition71111}. But by \cite[Theorem 7.1 (1)]{BM1} the maximal Hilbert-Samuel function of $X_{p}$ must occur at a distinguished point in some affine chart. \qedsymbol
\\
\\\textbf{Maximal strata.} Define the maximal strata of $H$
$$
S_H(X) := \bigcup_{\sigma \in \Sigma : H_{X, a_{\sigma}}=H} S_{X_{\sigma}}(a_{\sigma}),
$$
where we take the union over all distinguished points that have associated Hilbert-Samuel function $H$ and consider the associated fibered strata.

\begin{lemma}
    \label[lemma]{lemma7.12.1.2} 
    $S_H(X) \subset Sing_{\mathbb{Z}}(X)$.
\end{lemma}

\textit{Proof.} Follows from \cref{corollary71124}. \qedsymbol

\begin{lemma}
    $S_H(X)$ can be equipped with a closed $\Lambda$-subscheme structure of $X$ and $M$.
\end{lemma}

\textit{Proof.} We argue as in the hypersurface case in \cref{lemma7825}. For each $\sigma \in \Sigma$ with $U_{\sigma} \cap S_H(X) \not = \emptyset$ let $L_{\sigma}$ be the ideal defined as in \cref{corollary71126} by Hasse derivative monomials and for all other $\sigma$ define $L_{\sigma}$ as the unit ideal. We show these $L_{\sigma}$ glue along faces. Let $\sigma \in \Sigma$ such that $U_{\sigma} \cap S_H(X) \not = \emptyset$ and $\tau$ a face of $\sigma$ so $U_{\tau}$ is given by inverting in $U_{\sigma}$  each $x_i, i \not \in V_{\tau}$, i.e. $U_{\tau} \cong D(f) \subset U_{\sigma}$ for $f=\Pi_{i: i \not \in V_{\tau}} x_i$. Let $a_{\sigma}$ be the distinguished point of $X_{\sigma} \hookrightarrow U_{\sigma}$ and $a_{\tau}$ be the distinguished point of $X_{\tau} \subset U_{\tau}$. By \cite[Theorem 7.1 (1)]{BM1}, $U_{\tau} \cap S_H(X) \not = \emptyset \iff a_{\tau} \in S_{X_{\sigma}}(a_{\sigma})$. Suppose $a_{\tau} \not \in  S_{X_{\sigma}}(a_{\sigma})$ then there exists some $i$ with $1 \leq i \leq t$ such that either there exists a $x_j$ with nonzero power in $x^{\alpha^i}$ and $a_{\tau} \not \in V(x_j)$ or $ord_{a_{\tau}}x^{\beta^i} < |\alpha^i|$. In both cases the same arguments in \cref{lemma7825} show that $L_{\sigma}$ contains an element that becomes a unit in $D(f)$ and hence $(L_{\sigma})_f=(1)$ as desired. Now suppose $a_{\tau} \in  S_{X_{\sigma}}(a_{\sigma})$. Then for each $i=1, \dots, t$ the image of $F^i$ in $U_{\tau}$ under the localization will form part of the standard basis of $X_{\tau} \subset U_{\tau}$ with $|\alpha^i| > 1$, i.e. will define $X_{\tau \subset N_{\tau}}$. Again by the same arguments in \cref{lemma7825} the Hasse derivative monomials in $U_{\tau}$ will come from $L_{\sigma}$ and so $(L_{\sigma}) \cong L_{\tau}$. It follows that we can glue the $L_{\sigma}$ to an ideal sheaf that defines a closed subscheme of $M$. By \cref{corollary71126} this subscheme is the $\mathbb{Z}$-realization of an equivariant closed monoid subscheme of $X$ and the underlying set agrees with $S_H(X)$. \qedsymbol

\begin{lemma} 
    \label[lemma]{lemma71214}
    As sets
    $$
    S_H(X) = \bigsqcup_{(p) \in \textrm{Spec} \, \mathbb{Z}} S_{H_{p}}(X_{p})  \subset X,
    $$
    where $S_{H_{p}}(X_{p}):=\{b \in X_p: H_{X_p, b}=H_p   \} \subset X_p$ is the maximal Samuel stratum of $X_p$ (see \cite[Section 8]{BM1}).
\end{lemma}

\textit{Proof.} Suppose we have some $b \in X$ such that $H_{X_{p}, b_p}=H_{p}$, i.e. the base change of $b$ has maximal Hilbert-Samuel function and so $b_p \in S_{H_{p}}(X_{p})$. In an affine open $X_{\sigma}$ containing $b$ the corresponding distinguished point $a_{\sigma}$ satisfies (by \cite[Theorem 7.1]{BM1}) $H=H_{p}=H_{(X_{\sigma})_{p}, (a_{\sigma}, p)_p}=H_{X_{\sigma}, a_{\sigma}}$ and $b \in S_H(X)$. Conversely, by definition any $b \in S_H(X)$ in one of the charts containing it has the property that $H_{p}=H=H_{(X_{\sigma})_{p}, (b)_p}$. \qedsymbol
\\
\\A sequence of blow-ups as in \cref{theorem7201} is \textit{$H$-permissible} if conditions (1) and (2) are satisfied and in addition $C_j \subset S_H(X_j), j=0,\dots, t$, i.e. each $X_j$ has $\mathbb{Z}$-maximal Hilbert-Samuel function $H$ and each $C_j$ is contained in the maximal strata of $X_j$.


\begin{lemma} 
    \label[lemma]{lemma71215}
    A centre of blowing up, $D \subset M$, which satisfies (1) of the main theorem, is $H$-permissible if and only if for all prime ideals $(p)$, $(D)_p$ is $H_p$-permissible for the fiber $X_p \rightarrow M_p$ (see \cite[Section 8]{BM1}).
\end{lemma}

\textit{Proof.} By \cref{C7_7} and \cref{C7_3} the equivalence of conditions (1) and (2) to their analogues in \cite[Theorem 1.1]{BM1} is immediate. The equivalence of the maximal strata condition follows from \cref{lemma71214}. \qedsymbol

\begin{theorem} 
    (See \cite[Theorem 8.1]{BM1}).
    \label[theorem]{theorem71215} 
    There exists a $H$-permissible sequence of blow-ups such that $S_H(X_{t+1})=\emptyset$.
\end{theorem}

\textit{Proof of main \cref{theorem7201}.} Assume \cref{theorem71215}. By \cref{lemma71214}, $S_H(X_{t+1})=\emptyset$ implies $S_{H_{p}}((X_{t+1})_{p})=\emptyset$ for all $(p) \in \textrm{Spec} \, \mathbb{Z}$, which implies that the maximal Hilbert-Samuel function across each fiber has been reduced. By applying \cref{theorem71215} finitely many times and by the stabilization theorem applied to the maximal Hilbert-Samuel function of the fibers (see \cite[Theorem 8.1]{BM1}) this process terminates after finitely many steps when the transforms of each $X_{p}$ are regular at all closed points, hence smooth. It then follows from \cite[Proposition 6.5]{monoid} that the final transform of $X$ is smooth over $\textrm{Spec} \, \mathbb{Z}$. By \cref{lemma7.12.1.2}, condition (3) of \cref{theorem7201} is satisfied.  \qedsymbol


\subsection{Toric \texorpdfstring{$\Lambda$}{Λ}-marked monomial ideals}
In this section we prove \cref{theorem71215}, completing the proof of the main theorem. We define toric $\Lambda$-marked monomial ideals $\underline{\mathcal{H}}$ such that fibers over $(p)$ are marked monomial ideals of Bierstone and Milman \cite[Section 8.2]{BM1}. $\Lambda$-marked monomial ideals in \cite{OrderReductionLambdamarkedMonomialIdeals} were defined as a generalisation of the toric case and we use their order reduction in combination with \cite[Lemma 8.7]{BM1} to prove \cref{theorem71215}. 

\begin{definition}
    A \textit{toric $\Lambda$-marked monomial ideal} is a tuple
    $$
    \underline{\mathcal{H}}=(M, N, P, \mathcal{H}, e),
    $$
    where
    \begin{itemize}
        \item $M$ is a smooth toric $\Lambda$-scheme given by a regular fan $\Sigma_M$.
        \item $N$ is a smooth closed toric $\Lambda$-subscheme of $M$ induced by an unpointed closed embedding of monoid schemes. If $\Sigma_N$ is the corresponding regular fan, then $N \rightarrow M$ is induced by a morphism $\Sigma_N \rightarrow \Sigma_M$.
        \item $P$ is a smooth toric $\Lambda$-scheme given by an equivariant closed embedding of $N_{MSch}$. Equivalently, $P$ is a a single orbit closure of $N$ by \cref{C7_3}.
        \item $\mathcal{H}=\mathcal{H}_1+ \dots + \mathcal{H}_r \subset \mathcal{O}_M$, \textrm{where each }$\mathcal{H}_i \subset \mathcal{O}_M$ can be written as $\mathcal{H}_i = \Pi_j \mathcal{D}_{i,j}$, where each $\mathcal{D}_{i,j}$ defines a codimension one orbit closure of $M$ such that $N$ is transverse to $\mathcal{D}_{i,j}$ and $\mathcal{D}_{i,j}$ does not contain $P$ (see \cite[Section 3.4]{OrderReductionLambdamarkedMonomialIdeals} and \cite[Section 3.5]{OrderReductionLambdamarkedMonomialIdeals}).
        \item $e$ \textrm{ is a positive integer.}
    \end{itemize}
\end{definition}

\textbf{$\Lambda$-marked monomial ideals.} Let $\underline{\mathcal{H}}=(M, N, P, \mathcal{H}, e)$ be a toric $\Lambda$-marked monomial ideal and $E=\sum_i E^i$ be the locally toric simple normal crossings divisor of $M$ given by all codimension one orbit closures (see \cref{example7412}). By \cref{example7412}, \cref{remark7413}, \cref{example71001}, \cref{lemma71141} and the fibers $P_p, N_p, M_p$ being smooth toric varieties (hence connected)
$$
\underline{\mathcal{H}}_E:=(M, N, P, \mathcal{H}, E, e)
$$
is a $\Lambda$-marked monomial ideal (see \cite[Definition 4.1.1]{OrderReductionLambdamarkedMonomialIdeals}) and we can use the results of \cite[Section 4.2]{OrderReductionLambdamarkedMonomialIdeals}.

\begin{remark}
    In contrast to the marked monomial ideals of \cite[Definition 8.2]{BM1} we specify that for our toric $\Lambda$-marked monomial ideals $P$ is connected. We require this to use order reduction for $\Lambda$-marked monomial ideals. In \cite[Section 4.2]{OrderReductionLambdamarkedMonomialIdeals} this property ensures that on an affine open chart $U_{\sigma}$ of $M$ the variables corresponding to the codimension one orbit closures do not vanish on $P \cap U_{\sigma}$. Note that if $P$ was a disjoint union of orbit closures then a codimension one orbit closure could contain a component of $P$ while not containing the whole of $P$. 
\end{remark}

\textbf{Fibers of toric $\Lambda$-marked monomial ideals.} Let $(p) \subset \mathbb{Z}$ be a prime ideal and define the fiber of a toric $\Lambda$-marked monomial ideal $\underline{\mathcal{H}}$ to be the tuple $(\underline{\mathcal{H}}_E)_p$ but without the $E_p$ (see \cite[Definition 4.1.14]{OrderReductionLambdamarkedMonomialIdeals}), i.e.
$$
\underline{\mathcal{H}}_p=(M_p, N_p, P_p, \mathcal{H}_p, e).
$$

\begin{lemma}
    Let $\underline{\mathcal{H}}$ be a toric $\Lambda$-marked monomial ideal. Then for any $(p) \in \textrm{Spec} \, \mathbb{Z}$, the fiber $\underline{\mathcal{H}}_p$ is a marked monomial ideal in the sense of \cite[Definition 8.2]{BM1}.
\end{lemma}

\textit{Proof.} By \cref{C7_3}, $N_p \hookrightarrow M_p$ is a closed equivariant embedding of smooth toric varieties and $P_p$ is a smooth closed invariant subvariety of $N_p$ so it suffices to consider $\mathcal{H}_p$. But $E_p=\sum_i (E^i)_p$ will be a simple normal crossings divisor on $M_p$, where each $(E^i)_p$ corresponds to a codimension one orbit closure of $M_p$. The desired relations between $\mathcal{H}_p$, $P_p, N_p$ then hold by \cite[Definition 4.1.14]{OrderReductionLambdamarkedMonomialIdeals} and the fact that in fibers our transverse condition is equivalent to the simultaneous normal crossings condition in \cite[Definition 8.2]{BM1}. \qedsymbol

\begin{remark}
    As our toric $\Lambda$-marked monomial ideals $\underline{\mathcal{H}}$ have a connectedness assumption on $P$ required for resolution, the base changes $P_p$ are also connected. It follows that each marked monomial ideal $\underline{\mathcal{H}}_p$ will always have this additional assumption.
\end{remark}


\subsection{Transforms of toric \texorpdfstring{$\Lambda$}{Λ}-marked monomial ideals}

\begin{definition}
    (See \cite[Section 8.1]{BM1}).
    Let $\underline{\mathcal{H}}=(M, N, P, \mathcal{H}, e)$ be a toric $\Lambda$-marked monomial ideal. The \textit{support of $\underline{\mathcal{H}}$} is defined as 
    $
    supp\underline{\mathcal{H}} := supp\underline{\mathcal{H}}_E, 
    $
    the support of the corresponding $\Lambda$-marked monomial ideal
    (see \cite[Definition 4.1.2]{OrderReductionLambdamarkedMonomialIdeals}).
    It follows from \cite[Observation 4.1.15]{OrderReductionLambdamarkedMonomialIdeals} that as sets
    $$
    supp\underline{\mathcal{H}} = \bigsqcup_{(p) \in \textrm{Spec} \, \mathbb{Z}} supp\underline{\mathcal{H}}_p,
    $$
    where $supp\underline{\mathcal{H}}_p$ is as in \cite[Section 8.1]{BM1}. 
\end{definition}

\textbf{Transforms of toric $\Lambda$-marked monomial ideals.} Let $\underline{\mathcal{H}}=(M, N, P, \mathcal{H}, e)$ be a toric $\Lambda$-marked monomial ideal. Let $\pi_D:M' \rightarrow M$ be the blowing up with centre $D$  satisfying condition (1) of the main theorem, i.e. a disjoint union of $Z_{\Delta}$'s of $M$ with fibers given by the same $\Delta$'s defining orbit closures of $M_p$. $D$ will be a \textit{permissible centre} for $\underline{\mathcal{H}}$ and the induced $\pi$ will be called \textit{permissible} if in addition, $C:=D \cap N$ has smooth fibers and $C \subset supp \underline{\mathcal{H}}$. By \cref{lemma71141}, $C$ will be a disjoint union of orbit closures of $N$.

\begin{lemma} 
    \label[lemma]{lemma71232}
    A blow-up centre $D \subset M$ satisfying condition (1) is permissible for $\underline{\mathcal{H}}$ if and only if every fiber $D_{p}$ is permissible for $(\underline{\mathcal{H}})_{p}$ (see \cite[Section 8.1]{BM1}).
\end{lemma}

\textit{Proof.} Follows from the decomposition of $supp\underline{\mathcal{H}}$ into supports of fibers. \qedsymbol

\begin{remark}
    Let $\underline{\mathcal{H}}$ be a toric $\Lambda$-marked monomial ideal and $\underline{\mathcal{H}}_E$ the corresponding $\Lambda$-marked monomial ideal. A permissible blow-up centre for $\underline{\mathcal{H}}_E$ (see \cite[Definition 4.1.5]{OrderReductionLambdamarkedMonomialIdeals}) is a permissible centre for $\underline{\mathcal{H}}$. Moreover, each permissible centre for $\underline{\mathcal{H}}$ will be a disjoint union of permissible centres for $\underline{\mathcal{H}}_E$.
\end{remark}

\begin{definition}
    (See \cite[Definition 8.3]{BM1}).
    Let $\underline{\mathcal{H}}$ be a toric $\Lambda$-marked monomial ideal and $\pi_D: M' \rightarrow M$ be the blow-up by the permissible centre $D$. Let $\underline{\mathcal{H}}_E'=(M', N', P', \mathcal{H}',E', e')$ be the $\Lambda$-marked monomial ideal given by the sequence of permissible blow-ups defined by $D$ (see \cite[Definition 4.1.10]{OrderReductionLambdamarkedMonomialIdeals}). We define the \textit{transform} $\underline{\mathcal{H}}'$ of $\underline{\mathcal{H}}$ by $\pi$ as $\underline{\mathcal{H}}_E'$ without $E'$, i.e. the tuple $(M', N', P', \mathcal{H}', e')$. 
\end{definition}

\begin{lemma}
    \label[lemma]{lemma71235} 
    $\underline{\mathcal{H}}'$ is a toric $\Lambda$-marked monomial ideal. Moreover, for all prime ideals $(p) \subset \mathbb{Z}$, $(\underline{\mathcal{H}}')_p = (\underline{\mathcal{H}}_p)'$, i.e. taking the transform under a permissible blow-up commutes with taking the fiber then transforming as in \cite[Definition 8.3]{BM1}.
\end{lemma}

\textit{Proof.} It suffices to show the lemma when $D$ is given by a single orbit closure and $\underline{\mathcal{H}}_E'$ is the transform under a single blow-up. $M'$ is a blow-up given by star subdivision and hence corresponds to a smooth toric $\Lambda$-scheme. By \cref{lemma71141} $N'$ is a sequence of blow-ups given by star subdivisions and so is a smooth toric $\Lambda$-scheme. The closed embedding $N' \hookrightarrow M'$ is induced from an unpointed monoid scheme. It follows from the local description of $E'$ on blow-up charts $U_{\sigma'}$ given by cones $\sigma' \in \Sigma_{M'}$ that $E'$ will be the locally toric simple normal crossings divisor given by the codimension one orbit closures of $ M'$. Therefore, by the proof of \cite[Lemma 4.1.9]{OrderReductionLambdamarkedMonomialIdeals}, $P'$ is a single orbit closure of $N'$ and $\mathcal{H}'$ is an ideal given by codimension one orbit closures of $M'$ satisying the required conditions. Commuting with fibers follows from \cite[Lemma 4.1.17]{OrderReductionLambdamarkedMonomialIdeals}. \qedsymbol
\\
\\A \textit{permissible} sequence of blow-ups for $\underline{\mathcal{H}}$ will be a sequence of blow-ups
$$
M=M_0 \xleftarrow[]{\pi_1} M_1 \leftarrow \dots \xleftarrow[]{\pi_{t+1}} M_{t+1},
$$
where each $\pi_{j+1}, j = 0, \dots, t$ is a permissible blow-up for $\underline{\mathcal{H}}_j = (M_j, N_j, P_j, \mathcal{H}_j, e_j)$ and $\underline{\mathcal{H}}_{j+1} = (M_{j+1}, N_{j+1}, P_{j+1}, \mathcal{H}_{j+1}, e_{j+1})$ is the transform of $\underline{\mathcal{H}}_j$ by $\pi_{j+1}$, with $\underline{\mathcal{H}}_0 = \underline{\mathcal{H}}$. It follows from \cref{lemma71232} and \cref{lemma71235} that a sequence of blow-ups satisfying condition (1) of \cref{theorem7201} is permissible for $\underline{\mathcal{H}}$ if and only if for all prime ideals $(p)$ it induces a permissible sequence for the marked monomial ideal $\underline{\mathcal{H}}_p$ (see \cite[Definition 8.3]{BM1}).

\begin{definition}
    \label[definition]{definition71238} 
    A \textit{resolution of singularities} of $\underline{\mathcal{H}}$ is a sequence of permissible blow-ups such that $supp \underline{\mathcal{H}}_{t+1} = \emptyset$. It follows from the decomposition of $supp \underline{\mathcal{H}}$ into its fibers that a sequence of blow-ups satisfying condition (1) of \cref{theorem7201} is a resolution for $\underline{\mathcal{H}}$ if and only if it is a resolution for each 
    ${\underline{\mathcal{H}}}_p$ (\cite[Definition 8.3]{BM1}). Moreover, an order reduction for $\underline{\mathcal{H}}_E$ (see \cite[Definition 4.1.10]{OrderReductionLambdamarkedMonomialIdeals}) is a resolution for $\underline{\mathcal{H}}$ and vice versa.
\end{definition}


\subsection{Reduction to resolution of toric \texorpdfstring{$\Lambda$}{Λ}-marked monomial ideals}
Finally, we show how the existence of resolution of singularities for an arbitrary toric $\Lambda$-marked monomial ideal implies \cref{theorem71215} using \cite[lemma 8.7]{BM1}. Let $X \hookrightarrow M$ be the embedding of toric $\Lambda$-schemes in \cref{theorem7201} and let $H$ be the $\mathbb{Z}$-maximal Hilbert-Samuel function of $X$.
\\
\\\textit{Proof of \cref{theorem71215}.} From the order reduction in \cite[Section 4.2]{OrderReductionLambdamarkedMonomialIdeals} and \cref{definition71238} there exists resolution of singularities for toric $\Lambda$-marked monomial ideals. We have analogues of \cite[Example 8.4]{BM1} defined by standard basis elements and $N$ as in \cref{C7_11_2} and using \cref{lemma71215}. The arguments of \cite[lemma 8.7]{BM1} can then be applied. \qedsymbol

\begin{remark}
    (Relation to \cite[Theorem 14.1]{monoid}).
    Throughout each step in the proof of the resolution for $X \hookrightarrow M$ each blow-up centre induces blow-up centres for the fibers via base change. In particular, the way in which we have defined the analogous data; toric $\Lambda$-marked monomial ideals, $\mathbb{Z}$-maximal Hilbert-Samuel functions and inductive arguments imply that the fibers of the blow-up sequence are precisely a blow-up sequence for $X_p \hookrightarrow M_p$ given by applying the original Bierstone and Milman theorem \cite[Theorem 1.1]{BM1}. The blow-up sequence of monoid schemes in the proof of \cite[Theorem 14.1]{monoid} using $k$-realization ($k$ a field) and \cite[Theorem 1.1]{BM1} can then be achieved by first $\mathbb{Z}$-realizing (\cite[Remark 5.3.1]{monoid}), applying our resolution and then considering the underlying monoid schemes. More precisely, assuming $chark=p$, we have
    \begin{align*}
        (X_{MScheme})_k &= (X_{MScheme})_{\mathbb{Z}} \times_{\textrm{Spec} \, \mathbb{Z}} \textrm{Spec} \, k 
        \\&= (X_{MScheme})_{\mathbb{Z}} \times_{\textrm{Spec} \, \mathbb{Z}} \textrm{Spec} \, \mathbb{F}_p \times_{\textrm{Spec} \, \mathbb{F}_p} \textrm{Spec} \, k.
    \end{align*}
    
    It follows that the resolution given by applying Bierstone and Milman's theorem to $(X_{MScheme})_k \rightarrow (M_{Mscheme})_k$ is the base change of a fiber of our resolution and hence independent of the characteristic of $k$.
\end{remark}

\begin{remark}
    In the proof of \cref{theorem7201} there were choices made for the blow-up centres so the remark above is based on the same choices being made in applying Bierstone and Milman's theorem to the fibers. \cite[Section 9]{BM1} gives an algorithmic way of choosing the blow-up centres and it is likely we can similarly define an algorithmic way of choosing the blow-up centres in  \cref{theorem7201} such that the base change of the centres will be given by the algorithmic choices of Bierstone and Milman (cf. \cite[Remark 4.2.17]{OrderReductionLambdamarkedMonomialIdeals}).
\end{remark}

\medskip
\section*{Acknowledgements}
The author would like to thank Prof. Christian Haesemeyer for his supervision, generosity and invaluable guidance throughout his doctoral studies and the writing of this paper. The author would also like to acknowledge the generosity of the David Lachlan Hay Memorial Fund, who's support through the University of Melbourne Faculty of Science Postgraduate Writing-Up Award made the writing of this paper possible.

\goodbreak
\label{Bibliography}


\bibliography{Bibliography}

\end{document}